\documentclass[12pt]{article}
\usepackage{amsfonts,amssymb}
\usepackage{amsmath,amsfonts,amsthm,amssymb,color}
\usepackage{psfrag}
\usepackage{graphics}
\usepackage{graphicx}

\newcommand{\comments}[1]{}
\newcommand{\ashhalf}{\renewcommand{\arraystretch}{1.12}}
\newcommand{\ashalf}{\renewcommand{\arraystretch}{1.25}}
\newcommand{\asa}{\renewcommand{\arraystretch}{1.5}}
\newcommand{\asaa}{\renewcommand{\arraystretch}{1.75}}
\newcommand{\asb}{\renewcommand{\arraystretch}{2}}
\newcommand{\asc}{\renewcommand{\arraystretch}{2.25}}

\newcommand{\mm}[1]{}
\newtheorem{theorem}{Theorem}

\newtheorem{lemma}{Lemma}

\newcommand{\rref}[1]{(\ref{#1})}

 \date{December 14, 2014}

\title{
 Ergodic Theorem for Stabilization of a Hyperbolic PDE Inspired by
Age-Structured Chemostat \thanks{The work of Malisoff and Krstic was supported by   US National Science Foundation Grants 1408295 and 1408376, respectively.}}
\author{
Iasson Karafyllis \thanks{I. Karafyllis is with the Department  of Mathematics, National Technical University of Athens,
Heroon Polytechneiou 9, 15780 Athens, Greece. {\tt iasonkar@central.ntua.gr}.}  \and  Michael Malisoff \thanks{M. Malisoff is
with the Department  of Mathematics,  303
Lockett Hall, Louisiana
State University, Baton Rouge, LA
70803-4918, USA. {\tt malisoff@lsu.edu}. }  \and   Miroslav Krstic  \thanks{M. Krstic is with
the
Department  of Mechanical and Aerospace Engineering,
University of California, San Diego,
La Jolla, CA 92093-0411, USA. {\tt krstic@ucsd.edu}. }}

\usepackage{amsfonts,graphics}

\begin{document}
\maketitle
\begin{abstract}
 We study a feedback stabilization problem for   a first-order hyperbolic partial differential equation. The problem is inspired
 by the stabilization of
  equilibrium age profiles for an age-structured chemostat, using the dilution rate as the control.
 Two distinguishing features of the  problem are that (a)  the PDE has a multiplicative (instead of an additive) input and (b)  the state is fed back to the inlet boundary. We provide a sampled-data feedback  that  ensures stabilization under arbitrarily sparse sampling and that  satisfies input constraints. Our
   chemostat feedback  does not require measurement of the  age profile, nor does it require exact knowledge of the model.\smallskip
   
\end{abstract}



\newtheorem{itclaim}{Claim}
\newenvironment{claim}{\begin{itclaim}\rm}{\end{itclaim}}
\pagestyle{myheadings}
\thispagestyle{plain}

\section{Introduction}

Age-structured models have been used
in mathematical biology and mathematical demography  for a long time. Models of age-structured populations are either discrete (resulting in Leslie matrix models) or continuous  (which produce the  McKendrick-von Foerster equation); see \cite{BC01,C94}. Age-structured models have also been used in mathematical economics and  environmental engineering \cite{BHY11}.
  The study of continuous age-structured models has focused on two different research directions, namely, the study of the dynamics of age-structured models, and optimal control problems. Optimal control problems for age-structured models have been used in \cite{BHY11,S14}  and optimality conditions were derived in \cite{BH14,FTV03}.
  \mm{Stability properties for the dynamics have been well studied  \cite{RHM08,RGG12}.   }The ergodicity problem has been studied in many works (such as \cite{I88a,I88b}) and many results have been presented for one or multiple continuous age-structured population models; see  for instance \cite{TK06}.

  Since the McKendrick-von Foerster  equation (which is also called the   Lotka-von Foerster equation)
  is a first-order hyperbolic partial differential equation (PDE)\mm{ with a non-local boundary condition}, it is reasonable to expect   recent results on controlling   hyperbolic PDEs (such as \cite{BC11,CVKB13,DMVK13,KK14,KS08,VKC11}) to apply to age-structured models. However, much of this literature uses boundary controls. For   chemostat models  (where the dilution rate is the most commonly used control  \cite{GR06,MMH08,RGG12}),
 it is natural to consider stabilization problems where the dilution rate is used to stabilize
   a specific age profile. The dilution rate control enters into the PDE directly, not at the boundary.
   Chemostats form the foundation of much
    current research in
bioengineering, ecology, and population biology, and are important in  biotechnological  processes such as waste water treatment plants \cite{S11,SW95}. Stability properties for the dynamics of chemostats have been well studied; see   \cite{RHM08,RGG12} and the references therein.

   In this paper, we study a simplified age-structured chemostat  model without an equation for the substrate concentration, i.e., we consider the substrate concentration to be constant. This is justified in the important case where the inlet concentration of the substrate is used to control the substrate concentration  (or the substrate concentration is slowly varying).
   Two distinguishing features of the control problem we consider  are (a) that the PDE has a multiplicative (instead of an additive) input and (b) that the state is fed back to the inlet boundary; see the system dynamics \rref{GrindEQ_16_}-\rref{GrindEQ_17_} below.
   Moreover, we do not require that we know the   age profile of the microorganism. Only the
value of a linear functional of the state profile (i.e., an output) is known at certain times
(namely, the sampling times); see \rref{GrindEQ_52_} below. The controller determines the value of
the dilution rate, and its key feature is the use of the natural logarithm of the output
value.
See \rref{GrindEQ_59_} for our control formula.
   It is a sampled-data feedback  for stabilizing  an arbitrary positive valued age profile, i.e., only sampled measurements are required, and the control  is   constant between sampling times.\mm{ This contrasts with the existing chemostat results where no sampling was allowed.}
   The feedback is valued in a pre-specified bounded interval, to incorporate  input constraints.\mm{ This contrasts with   \cite{TK06}, where the dilution rate is constant and where the goal
is to understand limit cycles, instead of feedback stabilization.}

   Other key novel features of our work relative to the existing control literature for hyperbolic PDEs are that we achieve global exponential stabilization for all positive valued initial age distributions, with arbitrarily sparse sampling, and that we do not require exact model knowledge. Our key stability estimate
   is in terms of the sup norm of
the logarithmic deviation of the state profile from the equilibrium age profile; see  \rref{GrindEQ_78_}.
   The proof of our main result uses the strong ergodic theorem and    the connection between hyperbolic PDEs and integral delay equations (IDEs) from our prior work \cite{KK14}. To our knowledge, this is the first time that the ergodic theorem has been used to solve a control problem.
     Our simulations  show   good performance of our control under three operating conditions, and   so support our view that our work would be useful for industrial applications.

\subsection*{Definitions and Notation}
We use the following notation.
Let $\mathcal I\subseteq {\mathbb  R} $ be any interval and $\Omega \subseteq {\mathbb  R} ^{n} $ be any set. Let  $C^{0} (\mathcal I ; \Omega )$  be the class of all continuous functions $f:\mathcal I\to \Omega$, and  $C^{k} (\mathcal I ; \Omega )$ for any integer $k\ge 1$ be  the class of all  functions in $C^{0} (\mathcal I ; \Omega )$ all of whose partial derivatives up through order $k$ exist and are continuous on $\mathcal I$.
 Let $L^{\infty } (\mathcal I;\Omega )$
be the equivalence classes of all essentially bounded Lebesgue measurable functions
$f:\mathcal I\to \Omega$ with norm $||f||_\infty={\rm ess\  sup}_{a\in \mathcal I }|f(a)|$. Let $L^{1} (\mathcal I;\Omega )$  be the equivalence classes of measurable functions $f:\mathcal I\to \Omega$
 for which $\| f\| _{1}  <\infty $, where $\| f\| _{1} =\int_{\mathcal I}|f(t)|{\rm d}t$.
 For each   $x\in {\mathbb  R} $,  let $[x]$ be the integer part of $x$, i.e., the largest integer $p$ such that $p\le x$. We let $\mathcal K_\infty$ denote the set of all strictly increasing unbounded continuous functions $\kappa:[0,\infty)\to [0,\infty)$ such that $\kappa(0)=0$.

For any subset $S\subseteq \mathbb R$ and any $A> 0$, we let
$PC^1([0,A];S)$ denote the class of all continuous functions
$z:[0,A]\to S$ for which there exists a finite (or empty) subset $B$ subset of
$(0,A)$ such that:
(i) the derivative $({\rm d}z/{\rm d}a)(a)$  exists at every point in $(0,A)\setminus B$ and is a
continuous function on $(0,A)\setminus B$,
(ii) all right and left limits  of $({\rm d}z/{\rm d}a)(a)$ when $a$ tends to
a point  in the set $B
\cup\{0,A\}$ exist and are finite. Given any subset $S\subseteq \mathbb R$,
we let $PC^0(\mathcal I;S)$ denote the set of all piecewise continuous functions,
 i.e, the set of all
functions $u:\mathcal I\to \mathbb R$ for which there exists a (possibly empty) set $B\subseteq
\mathcal I$ such that:
(i) u is continuous on $\mathcal I\setminus B$,
(ii) the intersection of every bounded subset of $\mathcal I$ with $B$ is finite (or
empty), and
(iii) all right and left limits  of $u(t)$ when $t$ tends to
a point (from the right or from the left) in the set $B$
exist and are finite.

\section{Main Result for Controlled Age-Structured Model}\label{mr}

\subsection{Statement of Problem and Theorem}
We consider the age-structured chemostat model given by
\begin{equation} \frac{\partial  f}{\partial  t} (t,a)+\frac{\partial  f}{\partial  a} (t,a)=-\big(\mu (a)+D(t)\big)f(t,a)                        \label{GrindEQ_16_}\end{equation}
for all $(t,a)\in (0,\infty)\times (0,A)$ and
\begin{equation}
 f(t,0)=\int_{\scriptscriptstyle 0}^{\scriptscriptstyle A}k(a)f(t,a){\rm d}a  \; \; {\rm for\ all}\; t\ge 0\; ,  \label{GrindEQ_17_}
\end{equation}
 where    $A>0$ is any constant,  $\mu :[0,A]\to[0,\infty) $ and $k:[0,A]\to[0,\infty) $ are continuous functions, and  we assume that $\int_{\scriptscriptstyle 0}^{\scriptscriptstyle A}k(a){\rm d}a >0$. The system \rref{GrindEQ_16_}-\rref{GrindEQ_17_}
is a continuous age-structured model of a population in a chemostat.
the boundary condition \rref{GrindEQ_17_} is the renewal condition, which determines the number of newborn individuals $f(t,0)$ at each time $t\ge 0$,  $A>0$ is the maximum reproductive age,
 $\mu$ is   the mortality function,  $f$ is the density of the population of age $a\in [0,A]$ at time $t\ge 0$, and  $k$ is the birth modulus. Given any constants  $D_{\rm min} >0$ and $D_{\rm max} >D_{\rm min}$, the  variable $D\in PC^0 \left([0,\infty) ;[D_{\min } ,D_{\max } ]\right)$ is called the {\em dilution rate} and is the control. Physically meaningful solutions of \rref{GrindEQ_16_}-\rref{GrindEQ_17_} are those  satisfying $f(t,a)\ge 0$ for all $(t,a)\in[0,\infty) \times [0,A]$.

We assume that there is a constant $D^{*} \in (D_{\rm min},D_{\rm max})$ such that
\begin{equation} \label{GrindEQ_5_}
1=\int_{\scriptscriptstyle 0}^{\scriptscriptstyle A}\left\{k(a)\exp \left(-D^{*} a-\int_{\scriptscriptstyle 0}^a\mu (s){\rm d}s \right)\right\}{\rm d}a\; ,
\end{equation}
which  is the   Lotka-Sharpe equation \cite{BC01}.
 Our assumption that there is a constant $D^{*} \in (D_{\min } ,D_{\max } )$ satisfying \rref{GrindEQ_5_}  is necessary for the existence of a non-zero equilibrium point for  \rref{GrindEQ_16_}-\rref{GrindEQ_17_}. In fact, for any constant $M> 0$, any function of the form
\begin{equation} f^{*} (a)=M\exp \left(-D^{*} a-\int _{ 0}^{a}\mu (s){\rm d}s \right)  \label{GrindEQ_18_}\end{equation}
for all $a\in [0,A]$
 is an equilibrium point for \rref{GrindEQ_16_}-\rref{GrindEQ_17_}. Therefore, there are a continuum of equilibria. This implies that the dynamics \rref{GrindEQ_16_}-\rref{GrindEQ_17_} cannot be made open-loop asymptotically stable
 to an equilibrium
 with the constant control $D(t)=D^*$, which is another motivation for our globally exponentially stabilizing feedback control design; see Section \ref{remarks} for more discussion on the equilibria.

It is natural to try to design a dilution rate controller $D$ for \rref{GrindEQ_16_} that is a function of values of the densities of the newborn individuals, i.e., \rref{GrindEQ_17_}. However, such measurements may not   be easy to obtain in practice. On the other hand, it is often possible to find a continuous function $p:[0,A]\to [0,\infty )$ such that we can measure
\begin{equation}y(t)=\int _{\scriptscriptstyle 0}^{\scriptscriptstyle A}p(a)f(t,a){\rm d}a                                                     \label{GrindEQ_52_}
\end{equation}
at each time $t\ge 0$.
 For instance, the case $p(a)\equiv 1$ corresponds to measuring the concentration of the microorganisms. Given any desired positive constant lower and upper bounds $D_{\rm min}>0$ and $D_{\rm max}>D_{\rm min}$ for the controller, any constant $T>0$,  and any desired reference profile   \rref{GrindEQ_18_} for any   $M>0$, and setting
\begin{equation}\label{ystar}y^{*} =\int _{\scriptscriptstyle 0}^{\scriptscriptstyle A}p(a)f^{*} (a){\rm d}a,\end{equation}
we can prove that our age-structured chemostat dynamics
\rref{GrindEQ_16_}-\rref{GrindEQ_17_}, in closed loop with the piecewise defined control defined by
\begin{equation}D(t) =\max \left\{D_{\min } ,\min \left\{D_{\max } ,D^{*} +T^{-1} \ln \left(\frac{y(iT)}{y^{*} } \right)\right\}\right\}                    \label{GrindEQ_59_}\end{equation}
for all $t\in [iT,(i+1)T)$ and all integers $i\ge 0$, satisfies a uniform global asymptotic stability estimate
for all initial functions
$f_{\scriptscriptstyle 0} \in PC^1 \left([0,A];(0,\infty )\right)$  satisfying
 \begin{equation}\label{posd}f_{\scriptscriptstyle 0} (0)=\int _{\scriptscriptstyle 0}^{\scriptscriptstyle A}k(a)f_{\scriptscriptstyle 0} (a){\rm d}a, \end{equation}
which means that we require that $f(0,a)=f_{\scriptscriptstyle 0} (a)$ for all $a\in [0,A]$. Our main theorem is:

\begin{theorem}\label{th2}Let $A$, $T$,
$D_{\rm min}$, and $D_{\rm max}$ be any positive constants, with $D_{\rm max}>D_{\rm min}$.
 Let $\mu :[0,A]\to [0,\infty) $ and  $k:[0,A]\to [0,\infty) $ be any  continuous functions, and assume that   $k\in
 PC^1([0,A];[0,\infty))$ and   is not the zero function. Let $p:[0,A]\to [0,\infty)$ be any continuous function such that $\int _{0}^{\scriptscriptstyle A}p(a){\rm d}a >0$, and $D^*\in (D_{\rm min},D_{\rm max})$ be any constant satisfying the Lotka-Sharpe condition \rref{GrindEQ_5_}. Then there exist a constant $\sigma>0$ and a function $\kappa\in \mathcal K_\infty$ such that
for each function $f_{0} \in
PC^1([0,A];(0,\infty))$ satisfying \rref{posd}, the unique solution of \rref{GrindEQ_16_}-\rref{GrindEQ_17_} in closed loop with \rref{GrindEQ_59_}
with the initial condition
$f_0$
satisfies
\begin{equation}\max_{0\le a\le A} \big|\ln \left(f(t,a)/f^{*} (a)\right)\big|\; \le\;  \exp (-\sigma \, t)\, \kappa \left(\max_{0\le a\le A} \big|\ln \left(f_{0} (a)/f^{*} (a)\right)\big|\right)       \label{GrindEQ_78_}\end{equation}
for all $t\ge 0$.
\end{theorem}
\subsection{Discussion on Theorem \ref{th2}}\label{remarks}
We discuss the structure of the feedback control, as well as several key features that distinguish our controller analysis from existing results.

  The tracking error norm $\left|\ln \left(f(t,a)/f^{*} (a)\right)\right|$ in the statement of the theorem is motivated by the fact that \rref{GrindEQ_16_}-\rref{GrindEQ_17_} has the restricted state space
  $
  X=\{f \in PC^1([0,A];(0,\infty)): f(0)=\int_{\scriptscriptstyle 0}^{\scriptscriptstyle A} k(a)f(a){\rm d}a\}$. In fact, our logarithmic transformation
    $x(t,a)=\ln(f(t,a)/f^*(a))$ for $(t,a) \in [0,\infty)\times [0,A]$
  produces a control system whose state $x_t$ takes values in $\mathbb R$ and has equilibrium $x=0$,
   which is the usual setting for hyperbolic PDEs. When $M=1$, the state space in the new variables is  \begin{equation}\label{nonlinearity}\left\{x \in PC^1([0,A];\mathbb R):  {\rm exp}(x(0))=\int_{\scriptscriptstyle 0}^{\scriptscriptstyle A}
k(a)f^*(a){\rm exp}(x(a)){\rm d}a\right\}.\end{equation}  The nonlinear character of the control problem we consider here is illustrated by the nonlinearity exp in \rref{nonlinearity}.
\mm{lack of geometric and algebraic properties of the
above set (since it is neither a subspace of $PC^1([0,A];\mathbb R)$ nor a positive cone in
$PC^1([0,A];\mathbb R)$ nor a ball in $PC^1([0,A];\mathbb R)$, and  it is not convex).}

Consider the special case where the function $p$ in our output \rref{GrindEQ_52_} is the birth modulus $k$. Then
our output  is the  density  \rref{GrindEQ_17_} of  newborn individuals.
Hence, our   theorem is a general result on output feedback control that includes the special case where the density of newborn individuals is being measured. For our infinite dimensional systems, the state is the function $f_t$, i.e.,  $(f_{t})(a)=f(t,a)$ for all $a\in[0,A]$. Therefore, even if $p=k$, our output feedback is not a state feedback. The motivation for selecting a specific $M>0$ so that the desired
equilibrium is   \rref{GrindEQ_18_} is that for the complete chemostat (which also
has an equation for the substrate), the selected $M>0$ maximizes the yield,
in the context of  anaerobic digestion; see \cite[Section 2.4]{KKSL08}.

 While there is no explicit Lyapunov functional in our proof of Theorem \ref{th2}, the function $V_0(t)=|\ln(f(t,0)/f^*(0))|$ acts as a Lyapunov-Razumikhin functional. This can be seen by showing that
$v_0(t)=\ln(f(t,0)/f^*(0))$ satisfies a suitable IDE; see \rref{GrindEQ_3_}, \rref{GrindEQ_4_}, and \rref{newrelate}.
 The function $V_0$ is a Lyapunov-Razumikhin functional for an IDE, instead of  a Lyapunov functional for a PDE. Our theorem provides bounds for solutions of the PDE in terms of the history of $V_0$, because  \begin{equation}|\ln(f(t,a)/f^*(a))|\le V_0(t-a)+\int_{t-a}^t|D(s)-D^*|{\rm d}s\end{equation} for all $t\ge A$ and $a\in [0,A]$; see \rref{GrindEQ_130_} below.
  See also \cite{KK14}, which also uses functions of the form $V_0(t)=|\ln(f(t,0)/f^*(0))|$ and which then builds a Lyapunov-Krasovskii functional of the form $W(t)=\max \{ {\rm exp}(ps)V_0(t-s) : s\in[0,A] \}$ for a suitable $p\in \mathbb R$; see the proof of \cite[Theorem 2.6]{KK14}. In the original coordinates, this functional is $Q(t)=\max\{ {\rm exp}(pa)q(a)|\ln(f(t,a)/f^*(a))| : a\in[0,A] \}$ for a suitable function $q$.

\mm{Typically, stabilizing feedbacks are designed in   a Lyapunov function analysis. To see the rationale for our   control  \rref{GrindEQ_59_}, first notice that we can use the ergodic theorem to write  $v(t)=\ln(f(t,0)/f^*(0))$ as
$v(t)=g_1(t,y(iT),D(iT))+g_2(t,f_0)$ for all $t\in  [iT,(i+1)T)$ for suitable functions $g_1$ and $g_2$. Then our   control \rref{GrindEQ_59_} has the property that  $g_1((i+1)T,y(iT),D(iT))=0$, so the value of  $V((i+1)T)$ is minimized
when the effects of the $g_2$ term are disregarded.}

As noted in the introduction, to our knowledge Theorem \ref{th2} is the first application of an ergodic theorem to solve a feedback stabilization problem. The proof of Theorem \ref{th2} has several steps. First, we use the ergodic theorem
and a transformation to an IDE to obtain estimates for $\ln(y(iT)/y^*)$. The control
 is designed such that $y(t)\to y^*$ as $t\to\infty$. In the next step, we produce an estimate relating  $\ln(y(iT)/y^*)$ and $\ln(f(t,0)/f^*(0))$, which we use
in the third step to show that $f(t,0)\to f^*(0)$ as $t\to \infty$. In the next step, we use our IDE
transformation to relate $\ln(f(t,a)/f^*(a))$ for all $a\in [0,A]$ to $\ln(f(t,0)/f^*(0))$.
Finally, we show that $f(t,a)\to f^*(a)$ for all $a\in [0,A]$. Our proof requires several lemmas, which
we turn to next.

\section{Background: Uncontrolled Age-Structured Models}We review the needed background from \cite{I88a, I88b} on uncontrolled chemostats, and other material from \cite{P83}, which we use to prove Theorem \ref{th2} below.
Let $A>0$ be any constant, let $\mu :[0,A]\to[0,\infty) $ and $k:[0,A]\to[0,\infty) $ be any continuous functions, and assume that $\int_{\scriptscriptstyle 0}^{\scriptscriptstyle A}k(a){\rm d}a >0$. Consider the initial value  problem given by the two equations
\begin{equation}\label{GrindEQ_1_}
\frac{\partial  z}{\partial  t} (t,a)+\frac{\partial  z}{\partial  a} (t,a)=-\mu (a)z(t,a)
\end{equation}
for all $ (t,a)\in (0,\infty)\times (0,A)$ and
\begin{equation}
z(t,0)=\int_{\scriptscriptstyle 0}^{\scriptscriptstyle A}k(a)z(t,a){\rm d}a   \label{GrindEQ_2_}
\end{equation}for all $t\ge 0$,
  with
  an initial condition $z(0,a)=z_{ 0} (a)$ for all $a\in [0,A]$.
 System \rref{GrindEQ_1_}-\rref{GrindEQ_2_} is a continuous age-structured model of a population in a closed ecosystem with no control, where  $\mu$ is   the mortality function,  $z$ is the density of the population of age $a\in [0,A]$ at time $t\ge 0$, and  $k$ is the birth modulus.
  Physically meaningful solutions  are those satisfying    $z(t,a)\ge 0$ for all $(t,a)\in[0,\infty) \times [0,A]$.
The following existence and uniqueness result follows from \cite[Proposition 2.4]{I88a} and  \cite[Theorems 1.3-1.4]{P83}:
\begin{lemma}\label{lm1}
For each absolutely continuous function $z_{ 0} \in C^{0} \left([0,A];{\mathbb  R} \right)$ such that   $z_{ 0} (0)=\int_{\scriptscriptstyle 0}^{\scriptscriptstyle A}k(a)z_{ 0} (a){\rm d}a $, there is a unique  function $z:[0,\infty)\times [0,A]\to \mathbb R$ that satisfies: (a)  For each $t\ge 0$, the function $z_t$ defined by $\left(z_{t} \right)(a)=z(t,a)$ for $a\in [0,A]$ is in $L^{1} \left([0,A];{\mathbb  R} \right)$, (b) the function $\Phi:[0,\infty)\to  L^{1} \left([0,A];{\mathbb  R} \right)$ defined by $\Phi(t)=z_t$ is continuously differentiable, (c) for each $t\ge 0$, the function
 $z_{t} \in L^{1} \left([0,A];{\mathbb  R} \right)$ is  absolutely continuous and satisfies   $z_{t} (0)=\int_{\scriptscriptstyle 0}^{\scriptscriptstyle A}k(a)z_{t} (a){\rm d}a $, and (d) equation \rref{GrindEQ_1_} holds for almost all $t>0$ and  $a\in (0,A)$. Moreover, if $z_{ 0} (a)\ge  0$ for all $a\in [0,A]$, then $z(t,a)\ge 0$ holds for all $(t,a)\in[0,\infty) \times [0,A]$.\hfill$\square$\end{lemma}

 We refer to $z$ or the function $\Phi$ from Lemma \ref{lm1} as
  the solution of \rref{GrindEQ_1_}-\rref{GrindEQ_2_}. We also use:
  \begin{lemma}\label{lm2}
 If $k\in PC^{1} \left([0,A];[0,\infty) \right)$, then for every $z_{ 0} \in PC^{1} \left([0,A];{\mathbb  R} \right)$ satisfying \begin{equation}z_{ 0} (0)= \int_{\scriptscriptstyle 0}^{\scriptscriptstyle A}k(a)z_{ 0} (a){\rm d}a ,\label{icc}\end{equation}  the function $z:[0,\infty)\times [0,A]\to \mathbb R$ from Lemma \ref{lm1}    is   $C^1$  on \[\mathcal S=\left\{(t,a)\in(0,\infty) \times (0,A)\, :\, t-a\not\in B\cup\{0,A\}\right\},\] where $B$ is the finite (or empty) set where the
derivative of $z_0$ is not defined, and
it satisfies  \rref{GrindEQ_1_} on $\mathcal S$ and equation \rref{GrindEQ_2_} for all $t\ge 0$.
 Also,
 \begin{equation} z(t,a)={\rm exp}\left(- \int_{\scriptscriptstyle 0}^a\mu (s){\rm d}s\right)v(t-a)          \label{GrindEQ_3_}\end{equation}
 holds for all $(t,a)\in [0,\infty) \times [0,A]$,  where  $v\in C^{0} \left([-A,\infty );{\mathbb  R} \right)\bigcap C^{1} \left((0,\infty );{\mathbb  R} \right)$ solves    \begin{equation}v(t)=\int_{\scriptscriptstyle 0}^{\scriptscriptstyle A}k(a){\rm exp}\left(- \int_{\scriptscriptstyle 0}^a\mu (s){\rm d}s\right)v(t-a){\rm d}a  \label{GrindEQ_4_}  \end{equation}for all $t\ge 0$ for   the  initial condition
  $v(-a)=\exp \left(\int_{\scriptscriptstyle 0}^{ a}\mu (s){\rm d}s \right)z_{ 0} (a)$ for all  $a\in (0,A]$.
 \hfill$\square$\end{lemma}

In the context of Lemma \ref{lm2}, the function $z_t$ is of class $PC^1$ for
every $t\ge 0$.
  The solution of   \rref{GrindEQ_4_} is found by differentiating both sides of  \rref{GrindEQ_4_} with respect to $t$, then applying integration by parts on the interval $[0,A]$, and then solving
\begin{equation} \dot{v}(t)=\tilde{k}(0)v(t)-\tilde{k}(A)v(t  -  A)  +  \int_{\scriptscriptstyle 0}^{\scriptscriptstyle A}\frac{d\tilde{k}}{da} (a)v(t  -  a){\rm d}a,  \end{equation}
where $\tilde{k}(a)=k(a)\exp \left(-\int_{\scriptscriptstyle 0}^{a}\mu (s){\rm d}s \right)$.
 Recalling that $\int_{\scriptscriptstyle 0}^{\scriptscriptstyle  A}k(a){\rm d}a>0$,
we also define the continuous functional
$P:L^{1} \left([0,A];{\mathbb  R} \right)\to {\mathbb  R} $ by
\begin{equation} \label{GrindEQ_7_}
  P(z_0)   =  \frac{\int_{\scriptscriptstyle 0}^{\scriptscriptstyle A}z_{ 0} (a)\big(\int_{a}^{\scriptscriptstyle A}k(s)\exp \left(\int_{s}^a\left(\mu (l)+D^{*} \right){\rm d}l \right){\rm d}s\big){\rm d}a  }{\int_{\scriptscriptstyle 0}^{\scriptscriptstyle A}ak(a)\exp \left(-\int_{\scriptscriptstyle 0}^a\left(\mu (l)+D^{*} \right){\rm d}l \right){\rm d}a }
\end{equation}
where $D^*$ is from the Lotka-Sharpe condition \rref{GrindEQ_5_}.
Recall the following strong ergodicity result, which follows from \cite[Section 3]{I88b}:

 \begin{lemma}\label{lm3}
 There are constants $\varepsilon >0$ and $K\ge 1$ such that for every absolutely continuous function $z_{ 0} \in C^{0} \left([0,A];{\mathbb  R} \right)$ that satisfies \mm{$dz_{ 0}/da \in L^{1} \left([0,A];{\mathbb  R} \right)$ }and \rref{icc}, the  solution of \rref{GrindEQ_1_}-\rref{GrindEQ_2_} satisfies
\begin{equation}
\label{GrindEQ_6_}\ashalf\begin{array}{l}
  \int_{\scriptscriptstyle 0}^{\scriptscriptstyle A}{\rm exp}(J(a)) \big|z(t,a) - {\rm exp}\left(D^*(t-a)  - J(a) \right)P(z_0) \big|{\rm d}a \;  \le\\  K {\rm exp}((D^{*} -\varepsilon )t) \int_{\scriptscriptstyle 0}^{\scriptscriptstyle A}{\rm exp}(J(a)) \left|z_{ 0} (a)\right|{\rm d}a  \; \, {\rm  for\  all}\; \, t\ge 0, \end{array}   \end{equation}
          where $J(a)=\int_{\scriptscriptstyle 0}^{ a}\mu (s){\rm d}s$.\hfill$\square$\end{lemma}

Lemma \ref{lm3} follows
by choosing
$(S(t)z_0)(a)=z_{t}(a)=z(t,a)$ for all $a\in [0,A]$ as the semigroup  in \cite{I88b}.
If \textit{$k\in PC^{1} \left([0,A];[0,\infty) \right)$}, then for every $z_{ 0} \in PC^{1} \left([0,A];{\mathbb  R} \right)$ that satisfies \rref{icc}, we define
\begin{equation}\phi (t)=\exp \left(-D^{*} t\right)v(t)-P(z_0)\; \; {\rm  for\  all}\; t\ge -A, \label{GrindEQ_8_}\end{equation}
   where $v\in C^{0} \left([-A,\infty );{\mathbb  R} \right)$ is the solution of  \rref{GrindEQ_4_} with the initial condition  \[v(-a)=\exp \left(\int_{\scriptscriptstyle 0}^{a}\mu (s){\rm d}s \right)z_{ 0} (a)\] for all $a\in (0,A]$.
   Then $\phi\in C^0((-A,0);\mathbb R)$. Also, \rref{GrindEQ_3_} and \rref{GrindEQ_6_} give
\begin{equation} \label{GrindEQ_9_}
\; \; \; \; \; \; \int_{\scriptscriptstyle 0}^{\scriptscriptstyle A}\exp \left(-D^{*} a\right)\big|\phi (t-a)\big|{\rm d}a \le K\exp \left(-\varepsilon \, t\right)\int_{\scriptscriptstyle 0}^{\scriptscriptstyle A}\exp \left(\int_{\scriptscriptstyle 0}^a\mu (s){\rm d}s \right)\big|z_{ 0} (a)\big|{\rm d}a
\end{equation}
 for all $t\ge 0$. Therefore, by setting
  \begin{equation} C  =  \max_{0\le a\le A} k(a)\exp \left(-\int_{\scriptscriptstyle 0}^a\mu (s){\rm d}s \right),\label{Cchoice}\end{equation}
  it follows from \rref{GrindEQ_4_} and \rref{GrindEQ_5_} that
\begin{equation} \label{GrindEQ_10_}\asaa\begin{array}{rcl}
|\phi (t)|&=&\left\vert\displaystyle\int_{\scriptscriptstyle 0}^{\scriptscriptstyle A}k(a)\exp \left(-D^{*} a  -  \displaystyle\int_{\scriptscriptstyle 0}^a\mu (s){\rm d}s \right)\phi (t-a){\rm d}a\right\vert\\
&\le & C\displaystyle \int_{\scriptscriptstyle 0}^{\scriptscriptstyle A}\exp \left(-D^{*} a\right)\left|\phi (t-a)\right|{\rm d}a\\& \le&    KC\exp (-\varepsilon \, t)\displaystyle\int_{\scriptscriptstyle 0}^{\scriptscriptstyle A}\exp \left(\int_{\scriptscriptstyle 0}^a\mu (s){\rm d}s \right)\left|z_{ 0} (a)\right|{\rm d}a
\end{array}\end{equation}holds for all $t\ge 0$.

\section{Key Lemma}
 Our proof of Theorem \ref{th2} will also use the following key lemma, which follows from our
 recent results in \cite{KK14}:
\begin{lemma}\label{lm4}  Let $G\in C^{0} \left(\left[0,A\right];[0,\infty) \right)$, set  $L=\int_{\scriptscriptstyle 0}^{\scriptscriptstyle A}G(a){\rm d}a$, and  let $\Delta \in (0,A)$ be any constant such that \begin{equation}
\int_{\scriptscriptstyle 0}^{\Delta }G(a){\rm d}a <1. \end{equation} If $L>1$, then for each  $x_{ 0} \in L^{\infty } \left([-A,0);\mathbb R \right)$, the solution $x\in L_{loc}^{\infty } \left([-A,\infty );{\mathbb  R} \right)$ of
 \begin{equation}\label{newide}
x(t)=\int_{\scriptscriptstyle 0}^{\scriptscriptstyle A}G(a)x(t-a){\rm d}a \end{equation}
  with the initial condition $x(a)=x_{ 0} (a)$ for $a\in [-A,0)$  satisfies
\begin{equation} \label{GrindEQ_14_}\asa\begin{array}{rcl}
 \min\{a_1,a_1b^{1+t/h}\}
 &\le&     \mathop{\inf }\limits_{-A\le a<0} x(t+a)\\& \le&       \mathop{\sup }\limits_{-A\le a<0} x(t+a)\; \le \;
 \max\{a_2,a_2b^{1+t/h}\} \end{array}
\end{equation}for all $t\ge 0$,
where
\begin{equation}\begin{array}{l}a_1=\inf_{-A\le a<0} x_0(a), \;   a_2=\sup_{-A\le a<0} x_0(a)\; ,\\
c=\int_0^{\Delta}G(a){\rm d}a, \; b=\frac{L-c}{1-c}\; ,\end{array}
\end{equation}
and
$h=\min \{\Delta ,A-\Delta \}$.
 \hfill$\square$\end{lemma}

For the proof of Lemma \ref{lm4}, see Appendix \ref{app1}.
A  consequence of Lemma \ref{lm4} is that for every $z_{ 0} \in PC^{1} \left([0,A];{\mathbb  R} \right)$ satisfying $z_{ 0} (0)=\int_{\scriptscriptstyle 0}^{\scriptscriptstyle A}k(a)z_{ 0} (a){\rm d}a $  and $z_{ 0} (a)>0$ for all $a\in [0,A]$, the   solution of \rref{GrindEQ_1_}-\rref{GrindEQ_2_} satisfies $z(t,a)>0$ for all $(t,a)\in[0,\infty) \times [0,A]$. To see why, note that if \begin{equation}\int_{\scriptscriptstyle 0}^{\scriptscriptstyle A}k(a)\exp \left(-\int_{\scriptscriptstyle 0}^a\mu (s){\rm d}s \right){\rm d}a \ge 1,\end{equation} then we may apply  Lemma \ref{lm4}  to \rref{GrindEQ_4_}, by choosing $G(a)=k(a){\rm exp}(-\int_{\scriptscriptstyle 0}^a\mu(s){\rm d}s)$. If   \begin{equation}\int_{\scriptscriptstyle 0}^{\scriptscriptstyle A}k(a)\exp \left(-\int_{\scriptscriptstyle 0}^a\mu (s){\rm d}s \right){\rm d}a <1,\end{equation} then set $x(t)=\exp (pt)v(t)$ for all $t\ge -A$ for any constant $p>0$,
 where $v$ satisfies
 \rref{GrindEQ_4_}. Then  \begin{equation} x(t)=\int_{\scriptscriptstyle 0}^{\scriptscriptstyle A}k(a)\exp \left(pa-\int_{\scriptscriptstyle 0}^a\mu (s){\rm d}s \right)x(t-a){\rm d}a \; {\rm  for\ all}\; t\ge 0\end{equation}and $\int_{\scriptscriptstyle 0}^{\scriptscriptstyle A}k(a)\exp \left(pa-\int_{\scriptscriptstyle 0}^a\mu (s){\rm d}s \right){\rm d}a \ge 1$ when $p$ is large enough, since $\int_{\scriptscriptstyle 0}^{\scriptscriptstyle  A}k(a){\rm d}a>0$.

Also,
\begin{equation} \; \; \; \; \begin{array}{l} \displaystyle\min_{0\le a\le A} \left(\exp \left(D^{*} a+\int_{\scriptscriptstyle 0}^a\mu (s){\rm d}s \right)z_{ 0} (a)\right)\le \displaystyle\min_{-A\le a\le 0} \left(P(z_0) +\phi (t+a)\right) \\ \le \displaystyle\max_{-A\le a\le 0} \left(P(z_0) +\phi (t+a)\right)\; \; \le \; \;  \displaystyle\max_{0\le a\le A} \left(\exp \left(D^{*} a+\int_{\scriptscriptstyle 0}^a\mu (s){\rm d}s \right)z_{ 0} (a)\right) \; \end{array}\! \!
                                     \label{GrindEQ_15_}\end{equation}
     for all $t\ge 0$. The inequalities \rref{GrindEQ_15_} are obtained by using Lemma \ref{lm4}, with $L=1$, $x(t)$ from Lemma \ref{lm4} taken to be ${\rm exp}(-D^*t)v(t)$,  and $G(a)$ taken to be the integrand in curly braces in \rref{GrindEQ_5_}.
Moreover, using our choice \rref{GrindEQ_18_} of $f^*$, our formula \rref{GrindEQ_3_} for the solutions $z(t,a)$
of the uncontrolled chemostat \rref{GrindEQ_1_}-\rref{GrindEQ_2_}, our formula \rref{GrindEQ_8_} for $\phi(t)$,
\rref{GrindEQ_15_}, and the fact that the solution of the controlled chemostat dynamics \rref{GrindEQ_16_}-\rref{GrindEQ_17_} with
  $D(t) \equiv D^{*}$ and any initial condition $f (0, a) = f_0(a)$ with $f_0\in PC^1([0,A];\mathbb R)$ satisfies
  \begin{equation}f(t,a)={\rm exp}(-D^{*} t)z(t,a)\end{equation} for all $(t,a)\in [0,\infty)\times [0,A]$ where $z(t,a)$ is the solution of \rref{GrindEQ_1_}-\rref{GrindEQ_2_} for the initial condition $z(0,a)=f_0(a)$, we obtain the following inequalities for all $t\ge 0$:
  \[\ashalf\begin{array}{rcl}\min_{a\in [a,A]}(f_0(a)/f_*(a))& \le& \min_{a\in [a,A]}(f(t,a)/f_*(a))\\& \le&  \max_{a\in [a,A]}(f(t,a)/f_*(a))\;  \le\;  \min_{a\in [a,A]}(f_0(a)/f_*(a)).\end{array}\]
  The preceding inequalities show that every equilibrium profile  \rref{GrindEQ_18_} for every choice of the constant $M>0$ is stable. However, since every neighborhood of an
equilibrium profile (in the $L^1$ norm or in the sup norm) contains infinitely many equilibria, each equilibrium profile is stable but not asymptotically stable (neutral stability).

\section{Proof of Theorem \ref{th2}}
Existence and uniqueness of the solution of the closed-loop system
\rref{GrindEQ_16_}-\rref{GrindEQ_17_} with the control \rref{GrindEQ_59_}
  can be established by the method of steps, as follows. First notice that by
 Lemma \ref{lm1}, the solution $z(t,a)$ of \rref{GrindEQ_1_}-\rref{GrindEQ_2_}
   with the control \rref{GrindEQ_59_} and the
   initial condition
$z_0=f_0$ exists for all $t\ge 0$, and Lemmas \ref{lm2} and \ref{lm4}  guarantee that $z_{t}$ is of class
$PC^1([0,A];(0,\infty)$ for all $t\ge 0$. Assume that the solution of  \rref{GrindEQ_16_}-\rref{GrindEQ_17_}, in closed loop with
\rref{GrindEQ_59_}, is defined on $[0,iT]$ for some non-negative integer $i$ and that
$f_t\in PC^1([0,A];(0,\infty))$ for all $t\in [0,iT]$. Then $D(t)$ can be
defined uniquely by \rref{GrindEQ_59_} on $[iT,(i+1)T)$, and $D$ is of class
$PC^0([0,(i+1)T);[D_{\rm min},D_{\rm max}])$. Moreover,
  the solution $f$ of \rref{GrindEQ_16_}-\rref{GrindEQ_17_} with the control \rref{GrindEQ_59_} satisfies
 \begin{equation}\label{newrelate}
f(t,a)=\exp  \left(-\int _{0}^{t}D(l){\rm d}l  \right)z(t,a)\end{equation} for all $(t,a)\in [0,\infty) \times [0,A]$ wherever
the solution $f$ is defined. Hence, we are in a
position to  uniquely  define $f(t,a)$ on $[iT,(i+1)T]\times [0,A]$. Notice that
$f_t$ is of class $PC^1([0,A];(0,\infty))$ for all t in $[0,(i+1)T]$. We can
continue this process to conclude that the solution of \rref{GrindEQ_16_}-\rref{GrindEQ_17_} with the control \rref{GrindEQ_59_} is defined for all $t\ge 0$ and satisfies $f_t\in PC^1([0,A];(0,\infty))$
for all $t\ge 0$.

 Using the fact that the solution of \rref{GrindEQ_16_}-\rref{GrindEQ_17_} with \rref{GrindEQ_59_} satisfies \rref{newrelate}
   for all $(t,a)\in [0,\infty) \times [0,A]$, our choice
  \rref{GrindEQ_52_} of the output gives
\begin{equation}\asaa\begin{array}{rcl}y(t)\! \! &=&\! \! \exp \left(-\int _{0}^{t}D(l){\rm d}l \right)\int _{0}^{A}p(a)z(t,a){\rm d}a\\& =& \exp \left(-\int _{0}^{t}D(l){\rm d}l \right)\int _{0}^{A}p(a){\rm exp}\left(-\int_0^a\mu(s){\rm d}s\right)v(t-a){\rm d}a
\\
 \! \!  &=&\! \! \exp \left(D^{*} t-\int _{0}^{t}D(l){\rm d}l \right)\\&&\, \, \times \int _{0}^{A}p(a)\exp \left(-D^{*} a-\int _{0}^{a}\mu (s){\rm d}s \right)\left(P(f_{0} )+\phi (t-a)\right){\rm d}a
    \end{array}                 \! \!      \! \!       \label{GrindEQ_79_}\end{equation}
  for all $t\ge 0$, by our choices of $v$ and $\phi$ from
    \rref{GrindEQ_4_} and \rref{GrindEQ_8_}, and the relationship \rref{GrindEQ_3_} between $z$ and $v$.
  Using
  \rref{GrindEQ_18_} and \rref{ystar}, we conclude that
\begin{equation} \label{GrindEQ_81_}\asa\begin{array}{l}
y^{*} =\int _{0}^{A}p(a)f^{*} (a){\rm d}a =M\beta,\; \; {\rm where}\\
\beta =\int _{0}^{A}p(a)\exp \left(-D^{*} a-\int _{0}^{a}\mu (s){\rm d}s \right){\rm d}a\; .\end{array}
\end{equation}
  Combining \rref{GrindEQ_79_} and  \rref{GrindEQ_81_} gives the following for all   $i\ge 0$ and  $t\in [iT, (i+1)T)$:
\begin{equation} \label{GrindEQ_83_}\asb\; \; \; \; \; \begin{array}{l}
\ln \left(\frac{y(t)}{y^{*} } \right)\\=\ln \left(\frac{y(iT)}{y^{*} } \right)+D^{*} (t-iT)-\int _{iT}^{t}D(l){\rm d}l + \ln \left(\frac{P(f_{0} )+\int _{0}^{A}g(a)\phi (t-a){\rm d}a }{P(f_{0} )+\int _{0}^{A}g(a)\phi (iT-a){\rm d}a } \right)\\
 =\ln \left(\frac{y(iT)}{y^{*} } \right)-(D_{i} -D^{*} )(t-iT)+ \ln \left(\frac{P(f_{0} )+\int _{0}^{A}g(a)\phi (t-a){\rm d}a }{P(f_{0} )+\int _{0}^{A}g(a)\phi (iT-a){\rm d}a } \right),\end{array}
\end{equation}
where
\begin{equation}g(a)=\beta ^{-1} p(a)\exp \left(-D^{*} a-\int _{0}^{a}\mu (s){\rm d}s \right)                          \label{GrindEQ_84_}\end{equation}
 for all $a\in [0,A]$ and
  \begin{equation}D_{i} =\max \left\{D_{\min } ,\min \left\{D_{\max } ,D^{*} +T^{-1} \ln \left(\frac{y(iT)}{y^{*} } \right)\right\}\right\}              \label{GrindEQ_86_}\end{equation}
   for all integers $i\ge 0$,
 and where the second equality in \rref{GrindEQ_83_} followed from the sampling structure of our controller \rref{GrindEQ_59_}.
We now set
\begin{equation}\asa\begin{array}{l}
x(t)=\ln \left(\frac{y(t)}{y^{*} } \right), \;  x_{i} =\ln \left(\frac{y(iT)}{y^{*} } \right),\; {\rm and}\\  u_{i} =\ln \left(\frac{P(f_{0} )+\int _{0}^{A}g(a)\phi ((i+1)T-a){\rm d}a }{P(f_{0} )+\int _{0}^{A}g(a)\phi (iT-a){\rm d}a } \right) \end{array}\label{GrindEQ_87_}\end{equation}
for all integers $i\ge 0$. Choosing the positive constant
 \begin{equation}\label{deltachoice}\begin{array}{l}\delta =\frac{1}{2} \min \left\{(D_{\max } -D^{*} )T,(D^{*} -D_{\min } )T\right\},\end{array}\end{equation}
we use  the following claim:
\begin{claim}\label{claim1}
The inequality
\begin{equation} \label{GrindEQ_88_}
\left|x_{i+1} \right|\le \left|x_{i} \right|-\min \left\{\left|x_{i} \right|,2\delta \right\}+\left|u_{i} \right|
\end{equation}
holds for all integers $i\ge 0$.\hfill$\square$\end{claim}

For the proof of Claim \ref{claim1}, see Appendix \ref{proofofclaim1}. We also require the following two claims, which we also prove in the appendices:
\begin{claim}\label{claim2}
For all integers $i\ge 0$, the inequalities
\begin{equation} \label{GrindEQ_89_}
\asa \begin{array}{rcl} x_{i} &\ge& \min \left\{0,x_{0} +i\left(D^{*} -D_{\min } \right)T\right\}\\&&+{\displaystyle\min_{k=0,...,i}} \left(\ln \left(\frac{P(f_{0} )+\int _{0}^{A}g(a)\phi (iT-a){\rm d}a }{P(f_{0} )+\int _{0}^{A}g(a)\phi (kT-a){\rm d}a } \right)\right)\; \; {\rm and}\\  x_{i} &\le& \max \left\{0,x_{0} -i\left(D_{\max } -D^{*} \right)T\right\}\\&&+\, {\displaystyle\max_{k=0,...,i}} \left(\ln \left(\frac{P(f_{0} ) \int _{0}^{A}g(a)\phi (iT-a){\rm d}a }{P(f_{0} )+\int _{0}^{A}g(a)\phi (kT-a){\rm d}a } \right)\right) \end{array}\! \! \!
\end{equation}
are satisfied.\hfill$\square$\end{claim}

\begin{claim}\label{claim3}
The inequalities\begin{equation}\; \; \; \asaa\begin{array}{l} {\min \left\{0,x(0)\right\}+\min_{k=0,...,[t/T]} \left(\ln \left(\frac{P(f_0) +\int _{0}^{A}g(a)\phi (t-a){\rm d}a }{P(f_0) +\int _{0}^{A}g(a)\phi (kT-a){\rm d}a } \right)\right)} \\ {\le x(t)\le \max \left\{0,x(0)\right\}+\max_{k=0,...,[t/T]} \left(\ln \left(\frac{P(f_0) +\int _{0}^{A}g(a)\phi (t-a){\rm d}a }{P(f_0) +\int _{0}^{A}g(a)\phi (kT-a){\rm d}a } \right)\right)} \end{array} \label{90a} \end{equation}
and
\begin{equation}\left|x(t)\right|\le \left|x_{[t/T]} \right|+\left|\ln \left(\frac{P(f_{0} )+\int _{0}^{A}g(a)\phi (t-a){\rm d}a }{P(f_{0} )+\int _{0}^{A}g(a)\phi (\left[t/T\right]T-a){\rm d}a } \right)\right|\label{90b}\end{equation}
hold for all $t\ge 0$.\hfill$\square$\end{claim}

We can combine estimate \rref{90a} with our bounds \rref{GrindEQ_15_} on $P(z_0)+\phi(t+a)$, our choice \rref{GrindEQ_18_} of $f^*$, and our choices of $g$ and $\beta$ in  \rref{GrindEQ_81_} and \rref{GrindEQ_84_} (which imply that $\int _{0}^{\scriptscriptstyle A}g(a){\rm d}a =1$) to get
\begin{equation} \label{GrindEQ_92_}
\left|x(t)\right|\le \left|x_{0} \right|+\ln \left(\frac{\max_{0\le a\le A} \left(f_{0} (a)/f^{*} (a)\right)}{\min_{0\le a\le A} \left(f_{0} (a)/f^{*} (a)\right)} \right)
\end{equation}for all $t\ge 0$.
The proof of \rref{GrindEQ_92_} uses the fact that the upper and lower bounds in \rref{GrindEQ_15_} are independent of $t$.
By \rref{GrindEQ_10_} and  \rref{GrindEQ_15_}, we have
\begin{equation} \label{GrindEQ_35_}
\; \; \; \; \left|\phi (t)\right|\le K^*\exp (-\varepsilon  t)\int_{\scriptscriptstyle 0}^{\scriptscriptstyle A}f_{ 0} (a){\rm d}a,\; \, {\rm  where}\; \, K^*=KC{\rm exp}\left(\int_{\scriptscriptstyle 0}^{\scriptscriptstyle A}\mu(s){\rm d}s\right),\;
\end{equation}
and where\mm{\begin{equation}\ashhalf\begin{array}{rcl} \displaystyle\min_{0\le a\le A} \left\{\exp \left(D^{*} a+\int_{\scriptscriptstyle 0}^{a}\mu (s){\rm d}s \right)f_{ 0} (a)\right\}&\le& \displaystyle\min_{-A\le a\le 0} \left(P(f_0) +\phi (t+a)\right) \\ &\le& \displaystyle\max_{-A\le a\le 0} \left(P(f_0) +\phi (t+a)\right)\\&\le& \displaystyle\max_{0\le a\le A} \left(\exp \left(D^{*} a+\int_{\scriptscriptstyle 0}^{a}\mu (s){\rm d}s \right)f_{ 0} (a)\right)\end{array} \label{GrindEQ_15_}\end{equation}
  for all $t\ge 0$,}
  $K$ and $\varepsilon$ are from Lemma \ref{lm3} and $C$ was defined in \rref{Cchoice}.
    Let $j$ be the smallest integer in $[ [A/T]+1,\infty)$
        such that
\begin{equation} \label{GrindEQ_93_}
K^{*} \left\| f_{0} \right\| _{1} \exp (-\varepsilon \, (jT-A))\le \frac{\exp (\delta )-1}{\exp (\delta )+1} P(f_{0} ).
\end{equation}
where  $\delta$ is from \rref{deltachoice}.
We need the following claim, which we prove in Appendix \ref{proofofclaim4}:
\begin{claim}\label{claim4}
For all integers $i\ge j$, we have
\begin{equation}\label{claim4a} \left|u_{i} \right|\le \delta \; \; {\rm and}\; \; \left|u_{i} \right|\le \frac{K^{*} \left\| f_{0} \right\| _{1} (\exp (\delta )+1)\exp (\varepsilon A)}{P(f_{0} )} \exp (-\varepsilon \, iT).\end{equation}
 Also, \begin{equation}\label{claim4b}\asa\begin{array}{l} \left|\ln \left(\frac{P(f_{0} )+\int _{0}^{A}g(a)\phi (t-a){\rm d}a }{P(f_{0} )+\int _{0}^{A}g(a)\phi (iT-a){\rm d}a } \right)\right|\\\le \frac{K^{*} \left\| f_{0} \right\| _{1} \left(\exp (\delta )+1\right)\exp (\varepsilon A)}{P(f_{0} )} \exp \left(-\varepsilon \, iT\right)\end{array}\end{equation} holds for all integers $i\ge j$ and all $t\ge iT$.\hfill$\square$\end{claim}

We next show that
\begin{equation} \label{GrindEQ_99_}
\exp \left(\left|x_{i+1} \right|\right)-1\le \exp (-\delta )\left(\exp \left(\left|x_{i} \right|\right)-1\right)+\exp \left(\left|u_{i} \right|\right)-1
\end{equation}
holds for all $i\ge j$.  When $\left|x_{i} \right|\le 2\delta $, we can use  \rref{GrindEQ_88_} to get  $\left|x_{i+1} \right|\le \left|u_{i} \right|$, which   implies \rref{GrindEQ_99_}. On the other hand, when $\left|x_{i} \right|>2\delta $, we conclude from \rref{GrindEQ_88_} from Claim \ref{claim1} that $\left|x_{i+1} \right|\le \left|x_{i} \right|-2\delta +\left|u_{i} \right|$. The previous inequality, in conjunction with the fact that $\left|u_{i} \right|\le \delta $ for all $i\ge j$ (which follows from Claim \ref{claim4}) gives
\[\begin{array}{rcl} \exp (\left|x_{i+1} \right|)-1&\le& \exp \left(\left|x_{i} \right|-2\delta +\left|u_{i} \right|\right)-1\\&=&\mm{\exp \left(\left|u_{i} \right|\right)-1+\exp \left(\left|x_{i} \right|-2\delta +\left|u_{i} \right|\right)-\exp \left(\left|u_{i} \right|\right) \\ &=&\exp \left(\left|u_{i} \right|\right)-1+\exp \left(\left|u_{i} \right|\right)\left(\exp \left(\left|x_{i} \right|-2\delta \right)-1\right) \\ =}\exp \left(\left|u_{i} \right|\right)-1+\exp \left(\left|u_{i} \right|-2\delta \right)\left(\exp \left(\left|x_{i} \right|\right)-1+1-\exp (2\delta )\right) \\ &\le& \exp \left(\left|u_{i} \right|\right)-1+\exp \left(\left|u_{i} \right|-2\delta \right)\left(\exp \left(\left|x_{i} \right|\right)-1\right) \\ &\le& \exp \left(\left|u_{i} \right|\right)-1+\exp \left(-\delta \right)\left(\exp \left(\left|x_{i} \right|\right)-1\right). \end{array}\]
Hence, \rref{GrindEQ_99_} holds for all $i\ge j$.
Using \rref{GrindEQ_99_} and induction, it follows that
\begin{equation} \label{GrindEQ_100_}\begin{array}{rcl}
\exp \left(\left|x_{i} \right|\right)-1&\le& \exp \left(-\delta (i-j)\right)\left(\exp \left(\left|x_{j} \right|\right)-1\right)\\&&+\, \sum _{l=j}^{i-1}\exp \left(-\delta (i-1-l)\right)\left(\exp \left(\left|u_{l} \right|\right)-1\right)\end{array}
\end{equation}
holds for all integers $i> j$.
\mm{More specifically, inequality \rref{GrindEQ_100_} follows from the definition of the sequence $\zeta _{i} ;=\exp \left(\left|x_{i} \right|\right)-1$ and the fact that inequality \rref{GrindEQ_99_} gives $\zeta _{i+1} \le \exp (-\delta )\, \zeta _{i} +\exp \left(\left|u_{i} \right|\right)-1$ for all $i\ge j$. Using induction, we can prove the formula \begin{equation}\zeta _{i} \le \exp \left(-\delta (i-j)\right)\, \zeta _{j} +\sum _{l=j}^{i-1}\exp \left(-\delta (i-1-l)\right)\left(\exp \left(\left|u_{l} \right|\right)-1\right) \end{equation} for all $i>j$, which directly implies \rref{GrindEQ_100_} for all $i>j$.}

Using our upper bounds \rref{claim4a} on $|u_i|$  from Claim \ref{claim4}
and the fact that ${\rm exp}(p)-1\le p{\rm exp}(p)$ for all $p\ge 0$, we get   $\exp \left(\left|u_{i} \right|\right)-1\le \exp (\delta )\left|u_{i} \right|$ for all $i\ge j$, and also the following consequence of  \rref{GrindEQ_100_}   for all $i>j$:
\[
\begin{array}{l} {\exp \left(\left|x_{i} \right|\right)-1\le \exp \left(-\delta (i-j)\right)\left(\exp \left(\left|x_{j} \right|\right)-1\right)} \\ {+\frac{K^{*} \left\| f_{0} \right\| _{1} \exp (\varepsilon A)(\exp (\delta )+1)}{P(f_{0} )} \displaystyle\sum _{l=j}^{i-1}\exp \left(-\delta (i-2-l)\right)\exp (-\varepsilon \, lT) }\; . \end{array}\]
Since $x\le \exp (x)-1\le x\exp (x)$ holds for all $x\ge 0$, we conclude that  the following holds for all $i>j$:
 \begin{equation} \label{GrindEQ_102_}
\! \! \asaa\begin{array}{rcl} \left|x_{i} \right|&\le& \mm{\exp \left(-\delta (i-j)\right)\exp \left(\left|x_{j} \right|\right)\left|x_{j} \right|  +\frac{K^{*} \left\| f_{0} \right\| _{1} \exp (\varepsilon A+2\delta )(\exp (\delta )+1)}{P(f_{0} )} \displaystyle\sum _{l=j}^{i-1}\exp \left(-\delta (i-l)\right)\exp (-\varepsilon \, lT)  \\ &\le&} \exp \left(-\tilde{\delta }(i-j)\right)\exp \left(\left|x_{j} \right|\right)\left|x_{j} \right| \\&& +\frac{K^{*} \left\| f_{0} \right\| _{1} \exp (\varepsilon A+2\delta )(\exp (\delta )+1)}{P(f_{0} )} \displaystyle\sum _{l=j}^{i-1}\exp \left(-\tilde{\delta }(i-l)\right)\exp (-\varepsilon \, lT),  \end{array}
\end{equation}
  where
  \begin{equation}\label{tildedeltachoice}
  \tilde{\delta }=\min \{\delta ,\varepsilon T\}.\end{equation} Since $\tilde{\delta }\le \varepsilon T$, it follows that $\exp (-\tilde{\delta }(i-l))\exp (-\varepsilon \, lT)\le \exp (-\tilde{\delta }\, i)$ for all $l=j,...,i-1$ and thus   \rref{GrindEQ_102_} implies  the following inequality for all $i>j$:
\begin{equation} \label{GrindEQ_103_}\ashalf\begin{array}{rcl}
\left|x_{i} \right|& \le&  \exp \left(-\tilde{\delta }(i-j)\right)\exp \left(\left|x_{j} \right|\right)\left|x_{j} \right|\\&&+\frac{K^{*} \left\| f_{0} \right\| _{1} \exp (\varepsilon A+2\delta )(\exp (\delta )+1)}{P(f_{0} )} (i-j)\exp \left(-\tilde{\delta }\, i\right)\end{array}
\end{equation}
  Notice that \rref{GrindEQ_103_} holds for $i=j$ as well and consequently, \rref{GrindEQ_103_} holds for all $i\ge j$. Using \rref{GrindEQ_92_} and \rref{GrindEQ_103_} and the fact that $x(iT)=x_i$ for all integers $i\ge 0$, we obtain the following inequality for all $i\ge j$:
\begin{equation} \label{GrindEQ_104_}
\asa\begin{array}{rcl}\; \;  \left|x_{i} \right|&\le&\frac{K^{*} \left\| f_{0} \right\| _{1} (\exp (\delta )+1)}{P(f_{0} )} \exp (2\delta +\varepsilon A)(i-j)\exp \left(-\tilde{\delta }\, i\right)\\ &&+\, \exp \left(-\tilde{\delta }(i-j)\right)\left(\left|x_{0} \right|+\ln \left(\frac{\max_{0\le a\le A} \left(f_{0} (a)/f^{*} (a)\right)}{\min_{0\le a\le A} \left(f_{0} (a)/f^{*} (a)\right)} \right)\right)\\&&\, \; \;  \; \; \; \times \exp \left(\left|x_{0} \right|\right)\left(\frac{\max_{0\le a\le A} \left(f_{0} (a)/f^{*} (a)\right)}{\min_{0\le a\le A} \left(f_{0} (a)/f^{*} (a)\right)} \right)
  \end{array}
\end{equation}
  Since $j$ is the smallest integer in $[ [A/T]+1,\infty)$
  that satisfies \rref{GrindEQ_93_} it follows that either (i) $j=[A/T]+1$ or (ii) $j>[A/T]+1$ and \begin{equation}K^{*} \left\| f_{0} \right\| _{1} \exp \big(-\varepsilon \, ((j-1)T-A)\big)\; >\; \frac{\exp (\delta )-1}{\exp (\delta )+1} P(f_{0} ).\end{equation} In either case, we have
  \begin{equation} \label{GrindEQ_105_}\; \; \; \; \; \;
\max \left\{1,\frac{K^{*} \left\| f_{0} \right\| _{1} }{P(f_{0} )} \right\}\frac{\exp (\delta )+1}{\exp (\delta )-1} \exp (\varepsilon \, (A+T))\; \ge\;  \exp (\varepsilon \, jT)\; \ge \; {\rm exp}(\tilde\delta j),
\end{equation}
by our choice \rref{tildedeltachoice} of $\tilde\delta$.

  Using \rref{GrindEQ_104_}-\rref{GrindEQ_105_} combined with the fact that \begin{equation}(i-j)\exp (-\tilde{\delta }\, i/2)\le i\exp (-\tilde{\delta }\, i/2)\le 2 \exp (-1)/\tilde\delta\end{equation} for all integers $i\ge j\ge 0$ (which follows because $r{\rm exp}(-r)\le {\rm exp}(-1)$ for all $r\ge 0$), we get the following inequality for all $i\ge j$:
\begin{equation} \label{GrindEQ_107_}\; \; \; \; \; \;
\left|x_{i} \right|\le \max \left\{1,\frac{K^{*} \left\| f_{0} \right\| _{1} }{P(f_{0} )} \right\}G\bigg(S(x_{0} ,f_{0} )\ln \left(S(x_{0} ,f_{0} )\right)+1\bigg)\exp \left(-\tilde{\delta }\, i/2\right),
\end{equation}where
\begin{equation} \label{GrindEQ_108_}
   \begin{array}{l} G=\big(\exp (\delta )+1\big)\exp (\varepsilon A)\max \left\{\frac{\exp (\varepsilon T)}{\exp (\delta )-1} ,\frac{2}{\tilde{\delta }} \exp (2\delta -1),\, 1\right\}\;  \, {\rm and}\\  S(x_{0} ,f_{0} )=\frac{\max_{0\le a\le A} \left(f_{0} (a)/f^{*} (a)\right)}{\min_{0\le a\le A} (f_{0} (a)/f^{*} (a))}{\rm exp}(|x_{0}|). \end{array}
\end{equation}
Using   \rref{90b}, \rref{GrindEQ_107_}, \rref{GrindEQ_108_}, the conclusion \rref{claim4b} from Claim \ref{claim4}, and our choice \rref{tildedeltachoice} of $\tilde\delta$, we obtain the following inequality  for all $t\ge jT$:
\begin{equation} \; \; \; \; \; \; \; \; \; \; \label{GrindEQ_109_}
\left|x(t)\right|\le \max \left\{1,\frac{K^{*} \left\| f_{0} \right\| _{1} }{P(f_{0} )} \right\}G\bigg(S(x_{0} ,f_{0} )\ln \left(S(x_{0} ,f_{0} )\right)+2\bigg)\exp \left(-\frac{\tilde{\delta }}{2} \left[\frac{t}{T} \right]\right)
\end{equation}
  Using \rref{GrindEQ_92_}, \rref{GrindEQ_105_} and \rref{GrindEQ_108_}, we get the following for all $t\in [0,jT]$:
\begin{equation} \label{GrindEQ_110_}
\asb\begin{array}{rcl} \left|x(t)\right|&\le& \left|x_{0} \right|+\ln \left(\frac{\max_{0\le a\le A} \left(f_{0} (a)/f^{*} (a)\right)}{\min_{0\le a\le A} \left(f_{0} (a)/f^{*} (a)\right)} \right)=\ln \left(S(x_{0} ,f_{0} )\right)\\&\le& \exp \left(-\frac{\tilde{\delta }}{2} \left[\frac{t}{T} \right]\right)\exp \left(\frac{\tilde{\delta }}{2} j\right)\ln \left(S(x_{0} ,f_{0} )\right) \\ &\le& \max \left\{1,\frac{K^{*} \left\| f_{0} \right\| _{1} }{P(f_{0} )} \right\}G\exp \left(-\frac{\tilde{\delta }}{2} \left[\frac{t}{T} \right]\right)S(x_{0} ,f_{0} )\ln \left(S(x_{0} ,f_{0} )\right) \end{array}
\end{equation}
 Estimate \rref{GrindEQ_110_} shows that inequality \rref{GrindEQ_109_} holds for all $t\ge 0$.

Defining
\begin{equation}\label{sigmachoice}\begin{array}{l}
\sigma =\frac{\tilde{\delta }}{4T} \; \; {\rm and}\; \; \tilde{G}=G\exp \left(\frac{\tilde{\delta }}{2} \right)\end{array}\end{equation}
 and using the fact that $\left[\frac{t}{T} \right]\ge \frac{t}{T} -1$, we can use  \rref{GrindEQ_92_}, \rref{GrindEQ_108_}, and \rref{GrindEQ_110_} to obtain the following for all $t\ge 0$:
\[\begin{array}{l}
\left|x(t)\right|\le \min \left\{\max \left\{1,\frac{K^{*} \left\| f_{0} \right\| _{1} }{P(f_{0} )} \right\}\tilde{G}\mathcal N(x_0,f_0)\exp \left(-2\sigma \, t\right)\, ,\, \ln \left(S(x_{0} ,f_{0} )\right)\right\},\end{array}\! \]
where
$\mathcal N(x_0,f_0)=S(x_{0} ,f_{0} )\ln \left(S(x_{0} ,f_{0} )\right)+2$.
 It now follows directly from  the fact that $\min \{a,b\}\le \sqrt{ab} $ for all $a\ge 0$ and $b\ge 0$ that
\begin{equation}\hspace{3.5em}\begin{array}{l}\left|x(t)\right|\le \left(\max \left\{1,\! \frac{K^{*} \left\| f_{0} \right\| _{1} }{P(f_{0} )} \right\}\! \tilde{G}\mathcal N(x_0,f_0)\ln \left(S(x_{0} ,f_{0} )\right)\right)^{1/2} \! \exp \left(-\sigma  t\right)      \end{array} \label{GrindEQ_112_}\end{equation}
  for all $t\ge 0$. Using \rref{GrindEQ_52_} and \rref{ystar}, we get
\begin{equation} \label{GrindEQ_113_}\; \; \; \; \; \; \; \;\;
y^{*} \min_{0\le a\le A} \left(f_{0} (a)/f^{*} (a)\right)\le y(0)=\int _{0}^{A}p(a)f_{0} (a){\rm d}a \le y^{*} \max_{0\le a\le A} \left(f_{0} (a)/f^{*} (a)\right)\; .
\end{equation}
 Using definition \rref{GrindEQ_87_} and \rref{GrindEQ_113_} then gives
\begin{equation} \label{GrindEQ_114_}
\left|x_{0} \right|\le \ln \left(\max \left\{\max_{0\le a\le A} \left(f_{0} (a)/f^{*} (a)\right),\frac{1}{\min_{0\le a\le A} \left(f_{0} (a)/f^{*} (a)\right)} \right\}\right).
\end{equation}

Next, we define the functions
\begin{equation} \label{GrindEQ_115_}\; \; \; \; \; \; \; \; \; \begin{array}{l}
Q(f_{0} )=\\\frac{\max_{0\le a\le A} \left(\frac{f_{0} (a)}{f^{*} (a)}\right)}{\min_{0\le a\le A} \left(f_{0} (a)/f^{*} (a)\right)} \max \left\{\max_{0\le a\le A} \left(f_{0} (a)/f^{*} (a)\right),\frac{1}{\min_{0\le a\le A} \left(f_{0} (a)/f^{*} (a)\right)} \right\}\end{array}
\end{equation}and
\begin{equation} \label{GrindEQ_116_}\; \; \; \; \;
\begin{array}{l} R(f_{0} )=\left(\max \left\{1,\frac{K^{*} \left\| f_{0} \right\| _{1} }{P(f_{0} )} \right\}\tilde{G}\bigg(Q(f_{0} )\ln \left(Q(f_{0} )\right)+2\right)\ln \left(Q(f_{0} )\bigg)\right)^{1/2}\; .\end{array}
\end{equation}

 It follows from   \rref{GrindEQ_112_} and our formula for $S$ from  \rref{GrindEQ_108_}
 that the following  holds:
\begin{equation}\left|x(t)\right|\le R(f_{0} )\exp \left(-\sigma \, t\right)\; \; {\rm  for\ all}\; t\ge 0\; .                                 \label{GrindEQ_117_}\end{equation}
   Using our formula for our control $D(t)$, definition \rref{GrindEQ_87_}, and \rref{GrindEQ_117_}, we obtain
\begin{equation} \; \; \; \;\left|D(t)-D^{*} \right|\; \le\;  T^{-1} R(f_{0} )\exp \left(-\sigma T\left[t/T\right]\right) \; \le\;  T^{-1} R(f_{0} )\exp \left(-\sigma (t-T)\right) \label{GrindEQ_118_}\end{equation}
for all $t\ge 0$, because $t/T\ge [t/T]-1\ge (t/T)-1$ for all $t\ge 0$.
Also, our relationship \rref{GrindEQ_3_} between $v(t-a)$ and the classical solution, combined with our formula
 \rref{GrindEQ_8_} for $\phi(t)$ and our relationship \rref{newrelate} between
 $f(t,a)$ and the solution $z(t,a)$ for the corresponding uncontrolled dynamics
 give
 \begin{equation}\label{tuvu3a}\begin{array}{l}
 P(f_0)={\rm exp}(-D^*t)v(t)-\phi(t)\\={\rm exp}(-D^*t)z(t,0)-\mathfrak{}\phi(t)={\rm exp}(-D^*t){\rm exp}\left(\int_{\scriptscriptstyle 0}^tD(\ell){\rm d}\ell\right)f(t,0)-\phi(t)\end{array}
 \end{equation}
for all $t\ge 0$. Hence, our output $y(t)$ satisfies
\[
\asa\begin{array}{l}y(t)= \int_{\scriptscriptstyle 0}^{\scriptscriptstyle  A}p(a)f(t,a){\rm d}a\; =\; \int_{\scriptscriptstyle 0}^{\scriptscriptstyle  A}p(a){\rm exp}\left(-\int_{\scriptscriptstyle 0}^tD(\ell){\rm d}\ell\right)z(t,a){\rm d}a\\
= \int_{\scriptscriptstyle 0}^{\scriptscriptstyle  A}p(a){\rm exp}\left(-\int_{\scriptscriptstyle 0}^tD(\ell){\rm d}\ell-\int_{\scriptscriptstyle 0}^a\mu(s){\rm d}s\right)v(t-a){\rm d}a\\
= \int_{\scriptscriptstyle 0}^{\scriptscriptstyle  A}p(a){\rm exp}\left(-\int_{\scriptscriptstyle 0}^tD(\ell){\rm d}\ell-\int_{\scriptscriptstyle 0}^a\mu(s){\rm d}s\right)\left(\phi(t-a)+P(f_0)\right){\rm exp}(D^*(t-a)){\rm d}a\\
= \int_{\scriptscriptstyle 0}^{\scriptscriptstyle  A}p(a){\rm exp}\left(-\int_{\scriptscriptstyle 0}^tD(\ell){\rm d}\ell-\int_{\scriptscriptstyle 0}^a\mu(s){\rm d}s\right)\left(\phi(t-a)-\phi(t)\right){\rm exp}(D^*(t-a)){\rm d}a\\\; \; \; \; +\, \int_{\scriptscriptstyle 0}^{\scriptscriptstyle  A}p(a){\rm exp}\left(-\int_{\scriptscriptstyle 0}^tD(\ell){\rm d}\ell-\int_{\scriptscriptstyle 0}^a\mu(s){\rm d}s\right){\rm exp}\left(-D^*a+\int_{\scriptscriptstyle 0}^tD(\ell){\rm d}\ell\right)f(t,0){\rm d}a\\
=\beta \, f(t,0)+\beta \, \exp \left(D^{*} t-\int _{0}^{t}D(l){\rm d}l \right)\int _{0}^{A}g(a)\left(\phi (t-a)-\phi (t)\right){\rm d}a
\end{array}\]
for all $t\ge 0$, where we used the relationship \rref{newrelate} between $z(t,a)$  and $f(t,a)$,
our choice \rref{GrindEQ_84_} of $g$, the relationship \rref{GrindEQ_3_} between $z(t,a)$ and $v(t-a)$, our choice \rref{GrindEQ_8_} of $\phi$, our formula in \rref{GrindEQ_81_} for $\beta$,
and the formula \rref{tuvu3a} for $P(f_0)$.
Dividing through  by $\beta M$ gives
\begin{equation} \label{GrindEQ_120_}\begin{array}{l}
f(t,0)/M=\\y(t)/y^{*} -M^{-1} \, \exp \left(D^{*} t-\int _{0}^{t}D(l){\rm d}l \right)\int _{0}^{A}g(a)\left(\phi (t-a)-\phi (t)\right){\rm d}a \, .\end{array}
\end{equation}
Using \rref{GrindEQ_35_}, \rref{GrindEQ_118_} and the fact that $\int _{0}^{\scriptscriptstyle A}g(a){\rm d}a =1$, we get  this for all $t\ge A$:
\begin{equation} \label{GrindEQ_121_}
\ashalf\begin{array}{l} \left|M^{-1} \, \exp \left(D^{*} t-\int _{0}^{t}D(l){\rm d}l \right)\int _{0}^{A}g(a)\left(\phi (t-a)-\phi (t)\right){\rm d}a \right|\\\le\;  2M^{-1} K^{*} \exp \left(-\varepsilon (t-A)\right)\left\| f_{0} \right\| _{1} \, \exp \left(\frac{1}{\sigma T} R(f_{0} )\exp (\sigma T)\right) \\ \le\;  2K^{*} \exp \left(-\varepsilon (t-A)\right)\frac{\left\| f_{0} \right\| _{1} \, }{\left\| f^{*} \right\| _{\infty } } \exp \left(\frac{1}{\sigma T} R(f_{0} )\exp (\sigma T)\right), \end{array}
\end{equation}
since $||f^*||_\infty\le M$. Using \rref{GrindEQ_18_},  \rref{GrindEQ_15_} with $f_0=z_0$, and \rref{GrindEQ_118_} gives
\begin{equation} \label{GrindEQ_122_}
\asa\begin{array}{l} {\left|M^{-1} \, \exp \left(D^{*} t-\int _{0}^{t}D(l){\rm d}l \right)\int _{0}^{A}g(a)\left(\phi (t-a)-\phi (t)\right){\rm d}a \right|\, } \\ \le \exp \left(\frac{1}{\sigma T} R(f_{0} )\exp (\sigma T)\right)\\\, \, \; \; \; \; \; \; \times \left(\max_{0\le a\le A} \left(f_{0} (a)/f^{*} (a)\right)-\min_{0\le a\le A} \left(f_{0} (a)/f^{*} (a)\right)\right) \end{array}
\end{equation}
for all $t\ge 0$ (by adding and subtracting $P(z_0)$ in the integrand in \rref{GrindEQ_122_}).

  Combining \rref{GrindEQ_121_} and \rref{GrindEQ_122_}, we obtain the following for all $t\ge 0$:
\begin{equation}\label{neweq}\asa\begin{array}{l} \left|M^{-1} \, \exp \left(D^{*} t-\int _{0}^{t}D(l){\rm d}l \right)\int _{0}^{A}g(a)\left(\phi (t-a)-\phi (t)\right){\rm d}a \right|\,  \\ \le \exp \left(\frac{1}{\sigma T} R(f_{0} )\exp (\sigma T)-\varepsilon (t-A)\right)\\\; \; \; \; \; \; \times \max \left\{\displaystyle\max_{0\le a\le A} \left(f_{0} (a)/f^{*} (a)\right)-\displaystyle\min_{0\le a\le A} \left(f_{0} (a)/f^{*} (a)\right),2K^{*} \frac{\left\| f_{0} \right\| _{1} \, }{\left\| f^{*} \right\| _{\infty } } \right\} \end{array}\end{equation}
We next define the following two functions:
\begin{equation} \label{GrindEQ_124_}\ashalf\; \; \; \; \;
\begin{array}{l} V(f_{0} )=\exp \left(\frac{1}{\sigma T} R(f_{0} )\exp (\sigma T)\right)\left(H(f_{0} )\max \left\{H(f_{0} ),2K^{*} \frac{\left\| f_{0} \right\| _{1} \, }{\left\| f^{*} \right\| _{\infty } } \right\}\right)^{1/2} \\ {\rm and}\; \, H(f_{0} )=\max_{0\le a\le A} \left(f_{0} (a)/f^{*} (a)\right)-\min_{0\le a\le A} \left(f_{0} (a)/f^{*} (a)\right)\; . \end{array}
\end{equation}
  Combining \rref{GrindEQ_122_} and \rref{neweq},
    and using the fact that $\min \{a,b\}\le \sqrt{ab} $ for all $a\ge 0$ and $b\ge 0$, gives
\begin{equation} \label{GrindEQ_125_}\ashalf\begin{array}{l}
\left|M^{-1} \, \exp \left(D^{*} t-\int _{0}^{t}D(l){\rm d}l \right)\int _{0}^{A}g(a)\left(\phi (t-a)-\phi (t)\right){\rm d}a \right|\\ \le V(f_{0} )\exp \left(-\varepsilon (t-A)/2\right)\end{array}
\end{equation}
for all $t\ge 0$.

  Using \rref{GrindEQ_3_}, our formula \rref{GrindEQ_8_} for $\phi(t)$, our bounds  \rref{GrindEQ_15_} with $z_0=f_0$, and our relationship \rref{newrelate} between $f(t,a)$ and the solution $z(t,a)$ of the corresponding uncontrolled system, we obtain the following for all $t\ge 0$:
\begin{equation}\begin{array}{l} \exp \left(-\int _{0}^{t}(D(l)-D^{*} ){\rm d}l \right)\displaystyle\min_{0\le a\le A} \left(f_{0} (a)/f^{*} (a)\right)\le\\ \frac{f(t,0)}{M}\;  \le\;  \exp \left(-\int _{0}^{t}(D(l)-D^{*} ){\rm d}l \right)\displaystyle\max_{0\le a\le A} \left(f_{0} (a)/f^{*} (a)\right)\end{array}
           \label{GrindEQ_126_}\end{equation}
 Using \rref{GrindEQ_118_} and \rref{GrindEQ_126_}, we obtain the following for all $t\ge 0$:
\begin{equation} \label{GrindEQ_127_}\ashalf
\begin{array}{l} \exp \left(-\mathcal A(f_0)\right)\displaystyle\min_{0\le a\le A} \left(f_{0} (a)/f^{*} (a)\right)\; \le\; \frac{f(t,0)}{M} \; \le\\  \exp \left(\frac{1}{\sigma T} R(f_{0} )\exp (\sigma T)\right)\displaystyle\max_{0\le a\le A} \left(f_{0} (a)/f^{*} (a)\right) \end{array}
\end{equation}
where $\mathcal A(f_0)=\frac{1}{\sigma T} R(f_{0})\exp (\sigma T)$.
 Using definition \rref{GrindEQ_87_} and \rref{GrindEQ_117_} we get the following for all $t\ge 0$:
\begin{equation} \label{GrindEQ_128_}
\exp \left(-R(f_{0} )\exp \left(-\sigma \, t\right)\right)\le \frac{y(t)}{y^{*} } \le \exp \left(R(f_{0} )\exp \left(-\sigma \, t\right)\right)
\end{equation}
Combining \rref{GrindEQ_120_}, \rref{GrindEQ_125_}, \rref{GrindEQ_127_}, \rref{GrindEQ_128_} and the fact that $\sigma \le \varepsilon /2$ (which follows from our choice \rref{sigmachoice} and the fact that  $\tilde{\delta }\le \varepsilon T$), we obtain the following for all $t\ge 0$:
\begin{equation}\ashalf\! \!   \begin{array}{l} {\max \! \left\{\exp \left(-\mathcal A(f_0)\right)\displaystyle\min_{0\le a\le A} \! \left(f_{0} (a)/f^{*} (a)\right),\mathcal R^\sharp(t,f_0)\right\}} \\ \le \frac{f(t,0)}{M} \le \min \! \left\{\exp \left(\mathcal A(f_0)\right)\displaystyle\max_{0\le a\le A} \! \left(f_{0} (a)/f^{*} (a)\right),\mathcal R^\sharp(t,f_0)\right\}, \end{array}
      \label{GrindEQ_129_}\end{equation}
where $\mathcal R^\sharp(t,f_0)=\exp \left(R(f_{0} )\exp (-\sigma \, t)\right)\! +\! V(f_{0} )\exp \left(-\sigma (t\! -\! A)\right)$.

  Using the relationship
  \rref{GrindEQ_3_} between the function $v$
   and the uncontrolled solution $z(t,a)$   of \rref{GrindEQ_1_}-\rref{GrindEQ_2_} with the initial condition $z(0,a)=f_{0} (a)$ and \rref{newrelate}, we obtain $v(t)=z(t,0)$, and therefore:
\begin{equation}\asa\; \; \; \; \begin{array}{rcl}
f(t,a)&=& {\rm exp}\left(-\int_{\scriptscriptstyle 0}^tD(\ell){\rm d}\ell-\int_{\scriptscriptstyle 0}^a\mu(s){\rm d}s\right)v(t-a)\\
&=&  {\rm exp}\left(-\int_{\scriptscriptstyle 0}^tD(\ell){\rm d}\ell-\int_{\scriptscriptstyle 0}^a\mu(s){\rm d}s\right)z(t-a,0)\\&
=& {\rm exp}\left(-\int_{\scriptscriptstyle 0}^tD(\ell){\rm d}\ell-\int_{\scriptscriptstyle 0}^a\mu(s){\rm d}s\right)f(t-a,0){\rm exp}\left(\int_{\scriptscriptstyle 0}^{t-a}D(\ell){\rm d}\ell\right)
\end{array}\end{equation}
when $t\ge a$, and
\begin{equation}\asa\; \; \; \begin{array}{l}
f(t,a)
= {\rm exp}\left(-\int_{\scriptscriptstyle 0}^tD(\ell){\rm d}\ell-\int_{\scriptscriptstyle 0}^a\mu(s){\rm d}s\right){\rm exp}\left(\int_{\scriptscriptstyle 0}^{a-t}\mu(s){\rm d}s\right)f_0(a-t)\end{array}\end{equation}
when $t\in [0,a)$, since $v(t-a)={\rm exp}(\int_{\scriptscriptstyle 0}^{a-t}\mu(s){\rm d}s)f_0(a-t)$.
 Hence, we can use the formula \rref{GrindEQ_18_} for $f^*$ to obtain
\begin{equation} \frac{f(t,a)}{f^{*} (a)}=\exp \left(-\int _{t-a}^{t}(D(l)-D^{*} ){\rm d}l \right)\frac{f(t-a,0)}{M} \; \; \label{GrindEQ_130_}\end{equation}for all   $(t,a)\in [0,\infty) \times [0,A]$ such that $t\ge a$,
and \begin{equation}\frac{f(t,a)}{f^{*} (a)}=\exp \left(-\int _{0}^{t}(D(l)-D^{*} ){\rm d}l  \right)\frac{f_{0} (a-t)}{f^{*} (a-t)}
\label{GrindEQ_131_} \end{equation}
for all  $(t,a)\in [0,\infty) \times [0,A]$ such that  $t< a$.
We now define the functions
\begin{equation} \label{GrindEQ_132_}
\ashalf \begin{array}{l} B_{1} (t,f_{0} )=\max \left\{\exp \left(-C(t,f_{0} )\right)-\Xi (t,f_{0} ),\right.\\\; \; \; \; \; \; \; \; \; \; \; \; \; \; \; \; \; \; \; \; \left.
\exp \left(-\frac{1}{\sigma T} R(f_{0} )\exp (\sigma T)\right)\min_{0\le a\le A} \left(f_{0} (a)/f^{*} (a)\right)\right\}, \\ B_{2} (t,f_{0} )=\min \left\{\exp \left(C(t,f_{0} )\right)+\Xi (t,f_{0} ),\right.\\\; \; \; \; \; \; \; \; \; \; \; \; \; \; \; \; \; \; \; \; \left.
\exp \left(\frac{1}{\sigma T} R(f_{0} )\exp (\sigma T)\right)\max_{0\le a\le A} \left(f_{0} (a)/f^{*} (a)\right)\right\}, \\ C(t,f_{0} )=R(f_{0} )\exp (-\sigma \, (t-A)),\; \; {\rm and}\\ \Xi (t,f_{0} )=V(f_{0} )\exp \left(-\sigma (t-2A)\right) \; .\end{array}\! \!
\end{equation}
  Combining \rref{GrindEQ_129_}, \rref{GrindEQ_130_}, \rref{GrindEQ_131_}, and \rref{GrindEQ_118_} gives the following for all $t\ge 0$:
\begin{equation} \label{GrindEQ_133_}\ashalf
\begin{array}{l} \exp \left(-AT^{-1} \sigma ^{-1} \exp (-\sigma (t-T-A))R(f_{0} )\right)B_{1} (t,f_{0} )\\\le \min_{0\le a\le A} \left(f(t,a)/f^{*} (a)\right)  \le \max_{0\le a\le A} \left(f(t,a)/f^{*} (a)\right)\\\le \exp \left(AT^{-1} \sigma ^{-1} \exp (-\sigma (t-T-A))R(f_{0} )\right)B_{2} (t,f_{0} ) \end{array}
\end{equation}

We now set\begin{equation}\begin{array}{l} w(t)=\displaystyle\max_{0\le a\le A} \left|\ln \left(f(t,a)/f^{*} (a)\right)\right| ,\\
w_{1} (t)=\displaystyle\max_{0\le a\le A} \ln \left(f(t,a)/f^{*} (a)\right)\; \; {\rm and}\\ w_{2} (t)=\displaystyle\max_{0\le a\le A} \ln \left(f^{*} (a)/f(t,a)\right)\end{array}
                                 \label{GrindEQ_134_}\end{equation}
for all $t\ge 0$. Clearly, definition \rref{GrindEQ_134_} implies that
\begin{equation}\asa\begin{array}{l} w(t)=\displaystyle\max_{0\le a\le A} \left|\ln \left(f(t,a)/f^{*} (a)\right)\right|\\=\displaystyle\max_{0\le a\le A} \left\{\displaystyle\max \left\{\ln \left(f(t,a)/f^{*} (a)\right),\ln \left(f^{*} (a)/f(t,a)\right)\right\}\right\} \\ {=\displaystyle\max \left\{\displaystyle\max_{0\le a\le A} \ln \left(f(t,a)/f^{*} (a)\right),\displaystyle\max_{0\le a\le A} \ln \left(f^{*} (a)/f(t,a)\right)\right\}} \\ {=\ln \left(\displaystyle\max \left\{\displaystyle\max_{0\le a\le A} \left(f(t,a)/f^{*} (a)\right),\displaystyle\max_{0\le a\le A} \left(f^{*} (a)/f(t,a)\right)\right\}\right)} \end{array} \label{GrindEQ_135_}\end{equation}for all $t\ge 0$, from which we get
\begin{equation}\begin{array}{l}\displaystyle\max_{0\le a\le A} \left(f(t,a)/f^{*} (a)\right)\le \exp (w(t))\\ {\rm and}\; \; \displaystyle\min_{0\le a\le A} \left(f(t,a)/f^{*} (a)\right)\ge \exp (-w(t))  \end{array}\label{GrindEQ_136_}\end{equation}
for all $t\ge 0$.
Therefore, our definitions \rref{GrindEQ_7_},  \rref{GrindEQ_115_}, and  \rref{GrindEQ_116_} for $P$, $Q$, and $R$ give
\begin{equation} \label{GrindEQ_138_}\asa\; \; \; \; \; \begin{array}{l}
P(f_{0} )\ge \exp (-w(0))P(f^{*}),\; \; Q(f_{0} )\le \exp (3w(0)),\; \; R(f_{0} )\le \tilde{Q}(w(0)),\; \;
{\rm and}\\ \left\| f_{0} \right\| _{1} \le \exp (w(0))\left\| f^{*} \right\| _{1},\end{array}
\end{equation}
where
\begin{equation}\tilde{Q}(s)=3\exp (3s)\left(\tilde{G}\max \left\{1,\frac{K^{*} \left\| f^{*} \right\| _{1} }{P(f^{*} )} \right\}\right)^{1/2} \left(\left(s+1\right)s\right)^{1/2}           \label{GrindEQ_141_}\end{equation}
for all $s\ge 0$.

Also, our definition of $V(f_0)$ in  \rref{GrindEQ_124_} in conjunction with our bounds \rref{GrindEQ_136_} and
\rref{GrindEQ_138_} give
\begin{equation} \label{GrindEQ_142_}
V(f_{0} )\le \tilde{P}(w(0)),
\end{equation}
where \begin{equation}\; \; \; \; \begin{array}{l}\tilde{P}(s)=\\\exp \left(\frac{1}{\sigma T} \tilde{Q}(s)\exp (\sigma T)\right)\left(2\sinh (s)\max \left\{2\sinh (s),2K^{*} \frac{\left\| f^{*} \right\| _{1} \, }{\left\| f^{*} \right\| _{\infty } } \exp (s)\right\}\right)^{1/2} \end{array}\label{GrindEQ_143_}\end{equation}
for all $s\ge 0$. Also, our definitions of $w_1$ and $w_2$ in \rref{GrindEQ_134_}
 in conjunction with estimate \rref{GrindEQ_133_} and our bound on $R(f_0)$ in \rref{GrindEQ_138_}
 give the following for all $t\ge 0$:
\noindent
\begin{equation} \label{GrindEQ_144_}
\asa\begin{array}{l} {w_{2} (t)\le AT^{-1} \sigma ^{-1} \exp (-\sigma (t-T-A))\tilde{Q}(w(0))+\ln \left(B_{1}^{-1} (t,f_{0} )\right)}\; \; {\rm and} \\ {w_{1} (t)\le AT^{-1} \sigma ^{-1} \exp (-\sigma (t-T-A))\tilde{Q}(w(0))+\ln \left(B_{2} (t,f_{0} )\right)} \end{array}
\end{equation}
Definitions \rref{GrindEQ_132_} in conjunction with \rref{GrindEQ_136_}, our bounds on $R(f_0)$ and $V(f_0)$ in \rref{GrindEQ_138_} and  \rref{GrindEQ_142_} and the facts that $\ln (a+b)\le \ln (a)+a^{-1} b$ for all $a>0$ and $b>0$ and  $\min \{a,b\}\le \sqrt{ab} $ for all $a\ge 0$ and $b\ge 0$, imply that
\begin{equation} \label{GrindEQ_145_}
\asaa\; \; \; \; \; \begin{array}{l} \ln \left(B_{2} (t,f_{0} )\right)\le\\ \min \left\{\frac{1}{\sigma T} \exp (\sigma T)\tilde{Q}(w(0))+w(0),\ln \left(\exp \left(C(t,f_{0} )\right)+\Xi (t,f_{0} )\right)\right\} \\ {\le \min \left\{\frac{1}{\sigma T} \exp (\sigma T)\tilde{Q}(w(0))+w(0),C(t,f_{0} )+\Xi (t,f_{0} )\exp \left(-C(t,f_{0} )\right)\right\}} \\ {\le \min \left\{\frac{1}{\sigma T} \exp (\sigma T)\tilde{Q}(w(0))+w(0),\exp \left(-\sigma (t-2A)\right)\left(R(f_{0} )+V(f_{0} )\right)\right\}} \\ \le \exp \left(-\sigma (t-2A)/2\right)\\\; \; \; \; \times \sqrt{\left(\frac{1}{\sigma T} \exp (\sigma T)\tilde{Q}(w(0))+w(0)\right)\left(\tilde{P}(w(0))+\tilde{Q}(w(0))\right)}  \end{array}
\end{equation}
holds for all $t\ge 0$. Also, \rref{GrindEQ_132_} in conjunction with \rref{GrindEQ_136_} and our bound \rref{GrindEQ_138_} on $R(f_0)$
give
\begin{equation}
\ln \left(B_{1}^{-1} (t,f_{0} )\right)\le \frac{1}{\sigma T} \tilde{Q}(w(0))\exp (\sigma T)+w(0)\; \; {\rm  for\ all}\; t\ge 0                \label{GrindEQ_146_}\end{equation}and
 \begin{equation}
\ashalf\begin{array}{l}
\ln \left(B_{1}^{-1} (t,f_{0} )\right)\le -\ln \big(\exp \left(-C(t,f_{0} )\right)-\Xi (t,f_{0} )\big)\;
\end{array}                               \label{GrindEQ_147_}\end{equation}
for all  $t\ge 0$ such that $1>\exp \left(C(t,f_{0} )\right)\Xi (t,f_{0} )$.
Using the facts that $\ln (1+x)\le x$ and  $e^{x} -1\le e^{x} x$ hold for all $x\ge 0$ and \rref{GrindEQ_147_}, we obtain the following for all $t\ge 0$ that satisfy $1\ge 2\exp \left(C(t,f_{0} )\right)\Xi (t,f_{0} )$:
\begin{equation} \label{GrindEQ_148_}
\asa\begin{array}{rcl} \ln \left(B_{1}^{-1} (t,f_{0} )\right)&\le& \ln \left(\frac{1}{\exp \left(-C(t,f_{0} )\right)-\Xi (t,f_{0} )} \right)\\&=&\ln \left(1+\frac{1-\exp \left(-C(t,f_{0} )\right)+\Xi (t,f_{0} )}{\exp \left(-C(t,f_{0} )\right)-\Xi (t,f_{0} )} \right) \\ &\le& \frac{\exp \left(C(t,f_{0} )\right)-1+\Xi (t,f_{0} )\exp \left(C(t,f_{0} )\right)}{1-\Xi (t,f_{0} )\exp \left(C(t,f_{0} )\right)} \\&\le& 2\left(\exp \left(C(t,f_{0} )\right)-1+\Xi (t,f_{0} )\exp \left(C(t,f_{0} )\right)\right)\\ &\le& 2\left(C(t,f_{0} )+\Xi (t,f_{0} )\right)\exp \left(C(t,f_{0} )\right) \end{array}
\end{equation}

Next note that that by our definitions \rref{GrindEQ_132_} and
our bounds on $R(f_0)$ and $V(f_0)$ in \rref{GrindEQ_138_} and  \rref{GrindEQ_142_},   the inequality $1\ge 2\exp \left(C(t,f_{0} )\right)\Xi (t,f_{0} )$ holds if \[1\ge 2{\rm exp}(\tilde Q(w(0)){\rm exp}(2\sigma A))\tilde P(w(0)){\rm exp}(-\sigma(t-2A)),\] which holds if
\begin{equation}\label{suffic}
0\ge \ln(2{\rm exp}(\tilde Q(w(0)){\rm exp}(2\sigma A))+\ln(\tilde P(w(0))+1)-\sigma(t-2A).\end{equation}
On the other hand, \rref{suffic} holds if
\[
t \ge 2A+\sigma^{-1}\ln(\tilde P(w(0))+1)+\sigma^{-1}\ln(2{\rm exp}(\tilde Q(w(0)){\rm exp}(2\sigma A)).\]
 Consequently, we conclude from  \rref{GrindEQ_148_}, \rref{GrindEQ_132_}, our bound on $R(f_0)$ from \rref{GrindEQ_138_},
    and \rref{GrindEQ_142_} that
\begin{equation}\; \; \; \begin{array}{l}\ln \left(B_{1}^{-1} (t,f_{0} )\right)\le\\ 2\exp \left(\tilde{Q}(w(0))-\sigma (t-2A)\right)\left(\tilde{P}(w(0))+\tilde{Q}(w(0))\right)\; {\rm  for\ all}\; t\ge \tilde{T}(w(0)),         \label{GrindEQ_149_}\end{array}\end{equation}
where $\tilde T(s)=2A+\sigma^{-1}\ln(\tilde P(s)+1)+\sigma^{-1}\ln(2{\rm exp}(\tilde Q(s){\rm exp}(2\sigma A))$
for all $s\ge 0$.
Then \rref{GrindEQ_146_} and \rref{GrindEQ_149_} give:
\begin{equation}\ashalf\begin{array}{l}\ln \left(B_{1}^{-1} (t,f_{0} )\right)\le \exp \left(-\sigma \left(t-\tilde{T}(w(0))\right)\right)\\\; \; \; \times \max \left\{w(0)+\frac{1}{\sigma T} \tilde{Q}(w(0))\exp (\sigma T),2\tilde{P}(w(0))+2\tilde{Q}(w(0))\right\} \end{array}\label{GrindEQ_151_}
\end{equation}for all $t\ge 0$.
 Also, \rref{GrindEQ_135_} and our definitions  of $w_1$ and $w_2$
  in \rref{GrindEQ_134_}
  give  $w(t)=\max \left\{w_{1} (t),w_{2} (t)\right\}$ for all $t\ge 0$.
       Using \rref{GrindEQ_144_}, \rref{GrindEQ_145_} and \rref{GrindEQ_151_} and noting that
(a) the functions $\tilde{P}$ and $\tilde{Q}$ in \rref{GrindEQ_143_} and  \rref{GrindEQ_141_}, respectively, are of class $\mathcal K_\infty $ and (b)   the function $\tilde{T}$ is non-decreasing, we conclude that there is a function $\kappa \in \mathcal K_\infty $ such that $w(t)\le \exp \left(-\sigma \, t/2\right)\kappa \left(w(0)\right)$ for all $t\ge 0$. Therefore, the   theorem follows from our definition of $w(t)$ from \rref{GrindEQ_134_}.

\section{Simulations}
To demonstrate our  control designs from Theorem  \ref{th2}, we carried out three simulations. In each simulation, we took the horizon $A=2$, the constant mortality function $\mu (a)=\mu =0.1$, $D^*=1$,
and the  birth modulus
\begin{equation}\begin{array}{l}
k(a)=
\left\{\begin{array}{lcl}ag,& a\in [0,1]\\
(2-a)g,& a\in [1,2]\end{array}\right.,\; \; {\rm where}\\[1em]
g=\frac{(\mu +D^{*} )^{2} }{\left(1-{\rm exp}(-(\mu +D^{*}) ) \right)^{2} } =2.718728\; .
\end{array}\end{equation}
The constant $g$ is chosen such that
  the Lotka-Sharpe condition
\rref{GrindEQ_5_} holds with $D^*=1$. The output is
\begin{equation} \label{GrindEQ__4a_}\begin{array}{l}
y(t)=\int _{0}^{2}f(t,a){\rm d}a\end{array}
\end{equation}
which  is the total concentration of the microorganism in the chemostat.
Our objective is to stabilize the equilibrium profile
\begin{equation}f^{*} (a)=\exp \left(-(D^{*} +\mu )a\right),\; \; a\in [0,2]\; .                                     \label{GrindEQ__5a_}
\end{equation} The equilibrium value of the output is
\begin{equation}\begin{array}{l}
y^{*} =\int _{ 0}^{2}f^{*} (a){\rm d}a= \frac{1-{\rm exp}(-2(D^*+\mu))}{D^*+\mu}=  0.808361.\end{array}\end{equation}
We tested   the output feedback law
\begin{equation}\; \; \; \; \; \; \; \begin{array}{l}D(t) =\\\max \left\{D_{\min } ,\min \left\{D_{\max } ,D^{*} +T^{-1} \ln \left(f(iT,0)/f^{*} (0)\right)\right\}\right\},\; \; t\in [iT,(i+1)T)\end{array}
   \label{GrindEQ__7a_}\end{equation}
  and the output feedback law
\begin{equation}\; \; \; \begin{array}{l}D(t) =\\\max \left\{D_{\min } ,\min \left\{D_{\max } ,D^{*} +T^{-1} \ln \left(y(iT)/y^{*} \right)\right\}\right\}, \; \; t\in [iT,(i+1)T) \end{array}\label{GrindEQ__8a_}
                         \end{equation}
                         where $i=0,1,2,\ldots$,
and for both controllers we chose $T=0.4$, $D_{\min } =0.5$, and $D_{\max } =1.5$.

We took our initial conditions to have the form
\begin{equation}f_{0} (a)=b_{0} -b_{1} a+c\exp (-\theta \, a),\; \; a\in [0,2]\label{inin}\end{equation}
where   $b_{0}$, $c$, and $\theta$ are positive parameters that we specify below, and where
\begin{equation} \label{GrindEQ__11_a}
b_{1} =g^{-1} (g-1)b_{0} +c\theta ^{-2} \left(1-{\rm exp}(-\theta ) \right)^{2} -cg^{-1}
\end{equation}
 is chosen so that   $f_{0} (0)=\int _{0}^{2}k(a)f_{0} (a){\rm d}a $ holds.
We must also choose the parameters such that  $\min_{a\in [0,2]}  f_{0} (a) >0$. We generated the simulations using a  uniform grid of function values
 $f(ih,jh)$ for
   $j=0,1,...,50$ and $i\ge 0$,  where $h=0.04$. For $i=0$ we had $f(0,jh)=f_{0} (jh)$ for $j=1,...,50$, where $f_0$ is from \rref{inin}.

   We computed the integrals
   \begin{equation}\label{tocompute}\begin{array}{l}
   y(ih)=\int _{0}^{2}f(ih,a){\rm d}a \; \; {\rm and}\; \; f(ih,0)=\int _{0}^{2}k(a)f(ih,a){\rm d}a \end{array}\end{equation} numerically for each $i\ge 0$. Since we wanted the numerical integrator to   evaluate   the integrals \rref{tocompute}
   exactly
   for every $i\ge 0$ when $f(ih,a)=C\exp (\sigma a)$ for certain real constants $C$ and $\sigma$, we did not use a conventional numerical integration scheme, such as the  trapezoid rule or Simpson's rule. The reason we wanted to evaluate the integrals exactly
    when $f(ih,a) $ is an exponential function is that the equilibrium profile  \rref{GrindEQ__5a_} is an exponential function and we would like to avoid a steady-state error due to the error induced by the numerical integrator. To this end, we set
    \begin{equation}\begin{array}{l}
    {\mathcal L}(i,j,h)=\ln \left(f(ih,(j+1)h)\right)-\ln \left(f(ih,jh)\right)\\ {\rm and}\; \; \mathcal I(i)=\{j: f(ih,(j+1)h) = f(ih,jh)\},\end{array}\end{equation} and we
    used the   integration schemes
   \begin{equation}\label{13aa}
   \int _{jh}^{(j+1)h}f(ih,a){\rm d}a \approx I_{i} (j)=\left\{
   \ashhalf\begin{array}{lll} h\frac{f(ih,(j+1)h)-f(ih,jh)}{{\mathcal L}(i,j,h)},& j\not\in \mathcal I(i)\\  hf(ih,jh), & j\in \mathcal I(i) \end{array}\right.\end{equation}
for $j=2,3,...,49$ and $i\ge 0$, and
\begin{equation} \label{GrindEQ__14_} \begin{array}{l}
\int _{0}^{2h}f(ih,a){\rm d}a \approx I_{i} (1)=\\[.5em]\left\{\ashhalf\begin{array}{lll} h\frac{f(ih,2h)-f^{2} (ih,h)/f(ih,2h)}{\ln \left(f(ih,2h)\right)-\ln \left(f(ih,h)\right)},& {\rm if}   \;  f(ih,h)\ne f(ih,2h) \\  2hf(ih,h),& {\rm if}\;  f(ih,h)=f(ih,2h) \end{array}\right.,\end{array}
\end{equation}
and we set
\begin{equation}\asc  \label{GrindEQ__15_} \; \; \;  \; \; \;    \begin{array}{l}\; \; \int _{jh}^{(j+1)h}af(ih,a){\rm d}a \approx J_{i} (j)=\\  \left\{\ashalf\ \! \! \! \! \!  \begin{array}{lll}
 h^{2}\left( \frac{f(ih,(j+1)h)+j\left(f(ih,(j+1)h)-f(ih,jh)\right)}{{\mathcal L}(i,j,h)} -\frac{  \left(f(ih,(j+1)h)-f(ih,jh)\right)}{{\mathcal L}^2(i,j,h)}\right),&   j\not\in \mathcal I(i) \\  \frac{2j+1}{2} h^{2} f(ih,jh),& j\in \mathcal I(i) \end{array}\right.\end{array}\! \! \! \! \! \! \! \end{equation}
for $j=2,3,...,24$ and $i\ge 0$, and
\begin{equation} \label{GrindEQ__16_} \; \; \; \; \; \;  \begin{array}{l} \textstyle\int _{0}^{2h}af(ih,a){\rm d}a \approx J_{i} (1)=\\[.25em]
\left\{\asa\begin{array}{lll} \frac{2h^{2} f(ih,2h)}{\ln \left(f(ih,2h)\right)-\ln \left(f(ih,h)\right)}  -\frac{h^{2} \left(f(ih,2h)-\frac{f^{2} (ih,h)}{f(ih,2h)} \right)}{\left(\ln \left(f(ih,2h)\right)-\ln \left(f(ih,h)\right)\right)^{2}}, & {\rm if}\; f(ih,h)\ne f (ih,2h)\\ 2h^{2} f(ih,2h),& {\rm if}\; f(ih,h)= f(ih,2h) \end{array}\right.\!\end{array}
\end{equation}and
 \begin{equation}\asc\; \; \;   \begin{array}{l}\int _{jh}^{(j+1)h}(2-a)f(ih,a){\rm d}a \approx K_{i} (j)=\\[.25em]
 \left\{\asa\begin{array}{lll} -\frac{h^{2} f(ih,(j+1)h)}{{\mathcal L}(i,j,h)}  +\left(2-jh+\frac{h}{{\mathcal L}(i,j,h)} \right)\frac{hf(ih,(j+1)h)-hf(ih,jh)}{ {\mathcal L}(i,j,h)}, &  j\not\in \mathcal I(i)\\  \left(2-\frac{2j+1}{2} h\right)hf(ih,jh), & j\in \mathcal I(i) \end{array}\right.\end{array}\label{GrindEQ__17_}\! \! \! \! \end{equation}
for $j=25,26,...,49$ and $i\ge 0$.

The derivation of formulas \rref{13aa}-\rref{GrindEQ__17_} is based on the interpolation of  \begin{equation}\tilde{f}_{j} (a)=C_{j} \exp (\sigma _{j} a)\end{equation} through the points $\left(jh,f(ih,jh)\right)$ and $\left((j+1)h,f(ih,(j+1)h)\right)$ for $j=1,2,...,49$. More specifically, we obtain the following  for $j=1,2,...,49$:
\begin{equation}\ashalf
\begin{array}{l} {\sigma _{j} =h^{-1} \ln \left(f(ih,(j+1)h)/f(ih,jh)\right)} \; \; {\rm and}\\  {C_{j} =f(ih,jh)\left(f(ih,(j+1)h)/f(ih,jh)\right)^{-j} } \end{array}          \label{stars}\end{equation}
  Using this interpolation, the exact integration formulas are used. For example, for \begin{equation}\int _{jh}^{(j+1)h}af(ih,a){\rm d}a \; {\rm  for}\; j=2,3,...,24,\end{equation} we get the following when $\sigma _{j} =h^{-1} \ln \left(f(ih,(j+1)h)/f(ih,jh)\right)\ne 0$:
\[\asaa\begin{array}{l} {\int _{jh}^{(j+1)h}af(ih,a){\rm d}a \approx \int _{jh}^{(j+1)h}a\tilde{f}_{j} (a){\rm d}a =C_{j} \int _{jh}^{(j+1)h}a\exp (\sigma _{j} a){\rm d}a } \\ {=C_{j} \sigma _{j}^{-1} h\exp (\sigma _{j} jh)\left((j+1)\exp (\sigma _{j} h)-j\right)-C_{j} \sigma _{j}^{-2} \left(\exp (\sigma _{j} (j+1)h)-\exp (\sigma _{j} jh)\right)} \end{array}\]
On the other hand, when $\sigma _{j} =h^{-1} \ln \left(f(ih,(j+1)h)/f(ih,jh)\right)=0$, we get
\[\int _{jh}^{(j+1)h}af(ih,a){\rm d}a \approx \int _{jh}^{(j+1)h}a\tilde{f}_{j} (a){\rm d}a=C_{j} \int _{jh}^{(j+1)h}a{\rm d}a  = \frac{h^{2} }{2} C_{j} \left(2j+1\right).\]
  Combining the above formulas with the estimated values for $C_{j}$ and $\sigma _{j} $ in \rref{stars}, we obtain formula \rref{GrindEQ__15_}. Similarly, we derive formulas \rref{13aa}, \rref{GrindEQ__14_}, \rref{GrindEQ__16_} and \rref{GrindEQ__17_}. Notice that the formulas \rref{13aa}, \rref{GrindEQ__14_}, \rref{GrindEQ__15_}, \rref{GrindEQ__16_}, and \rref{GrindEQ__17_} allow the numerical evaluation of the integrals \rref{tocompute} for every $i\ge 0$ without knowing $f(ih,0)$.

  Since the time and space discretization steps are both $h$,  we have   the exact formula
   \begin{equation}\ashalf\begin{array}{l}
   f((i+1)h,jh)=\\f(ih,(j-1)h)\exp \left(-(\mu +D_{i} )h\right)\; {\rm  for}\; j=1,2,...,50\; {\rm and\ all}\; i\ge 0 \; ,         \end{array}\label{GrindEQ__18_}\end{equation}
   where $D_i=D(ih)$.
Therefore, we have the following algorithm for simulating the closed-loop system:\medskip

  \underbar{Algorithm:} Given $f(ih,jh)$, for $j=1,...,50$ and certain $i\ge 0$ do the following:
\begin{enumerate}
 \item Calculate $f(ih,0)\approx g\sum _{j=2}^{24}J_{i} (j) +g\sum _{j=25}^{49}K_{i} (j) $, where $J_{i} (j)$ and $K_{i} (j)$ are given by \rref{GrindEQ__15_}, \rref{GrindEQ__16_}, and \rref{GrindEQ__17_}.

\item Calculate $y(ih)\approx \sum _{j=2}^{49}I_{i} (j) $, where $I_{i} (j)$ is given by \rref{13aa} and \rref{GrindEQ__14_}.

\item If $\frac{ih}{T} $ is an integer, then set \[D_{i} =\max \left\{D_{\min } ,\min \left\{D_{\max } ,D^{*} +T^{-1} \ln \left(y(ih)/y^{*} \right)\right\}\right\};\] otherwise, set $D_{i} =D_{i-1} $.

\item  Calculate $f((i+1)h,jh)$, for $j=1,...,50$ using \rref{GrindEQ__18_}.

\end{enumerate} \medskip

The above algorithm with obvious modifications was also used for the simulation of the open-loop system, and for the simulation of the closed-loop system under  the output feedback law \rref{GrindEQ__8a_}. We next present the results of our three simulations.

In our first simulation, we used the parameter values $b_{0} =0.2$, $b_{1} =0.15184212$, $c=0.8$, and $\theta =1$ in our initial conditions. In Figure 1, we plot the control values and the newborn individual values. We show the values for the open loop feedback $D(t)\equiv 1$, and for the state and output feedbacks from \rref{GrindEQ__7a_}
and \rref{GrindEQ__8a_}. Our simulation shows the efficacy of our control design.

\begin{figure}[t]\vspace{-4.5em}
\includegraphics[scale=0.5]{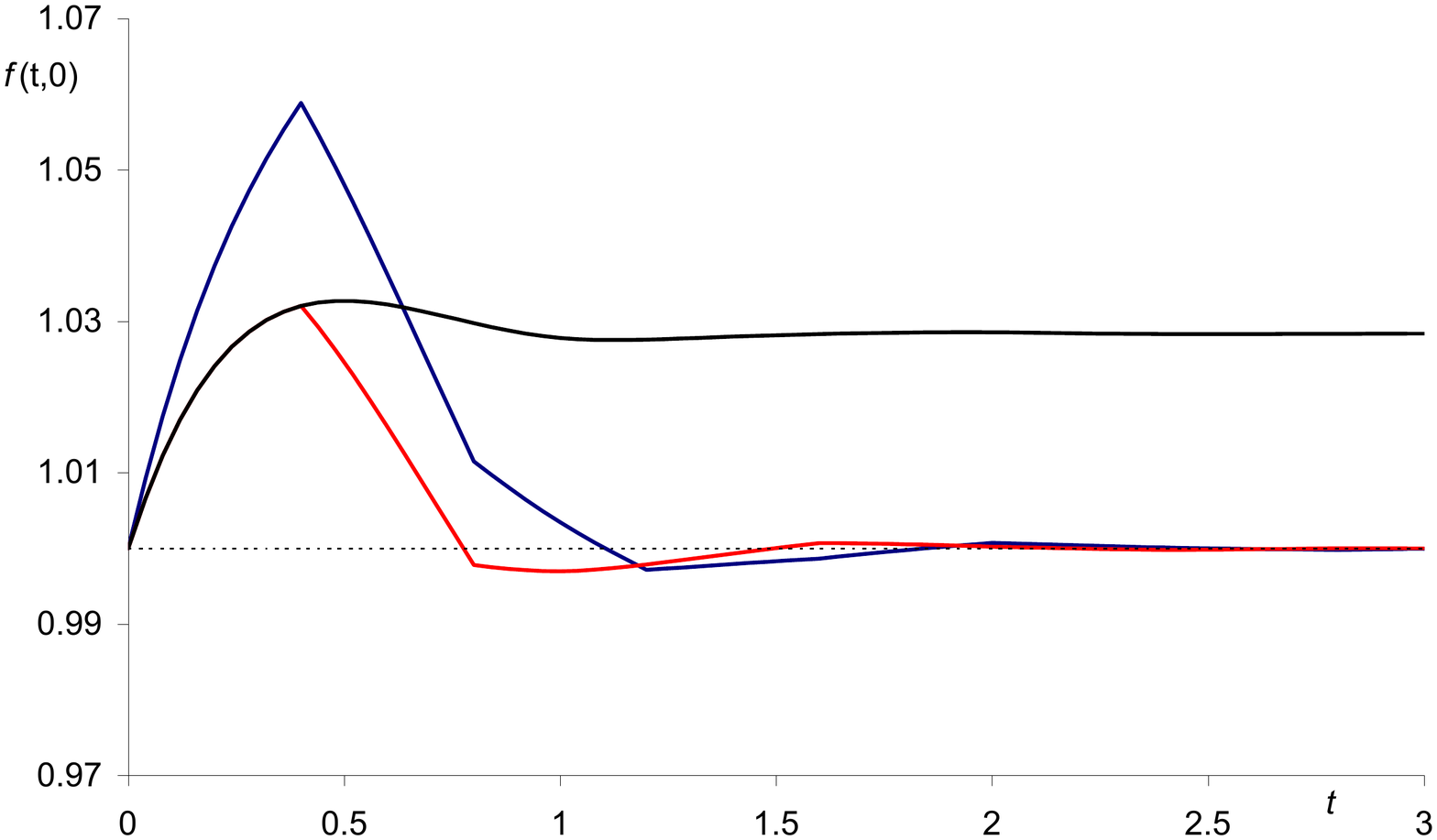}\\[-7em]
\includegraphics[scale=0.5]{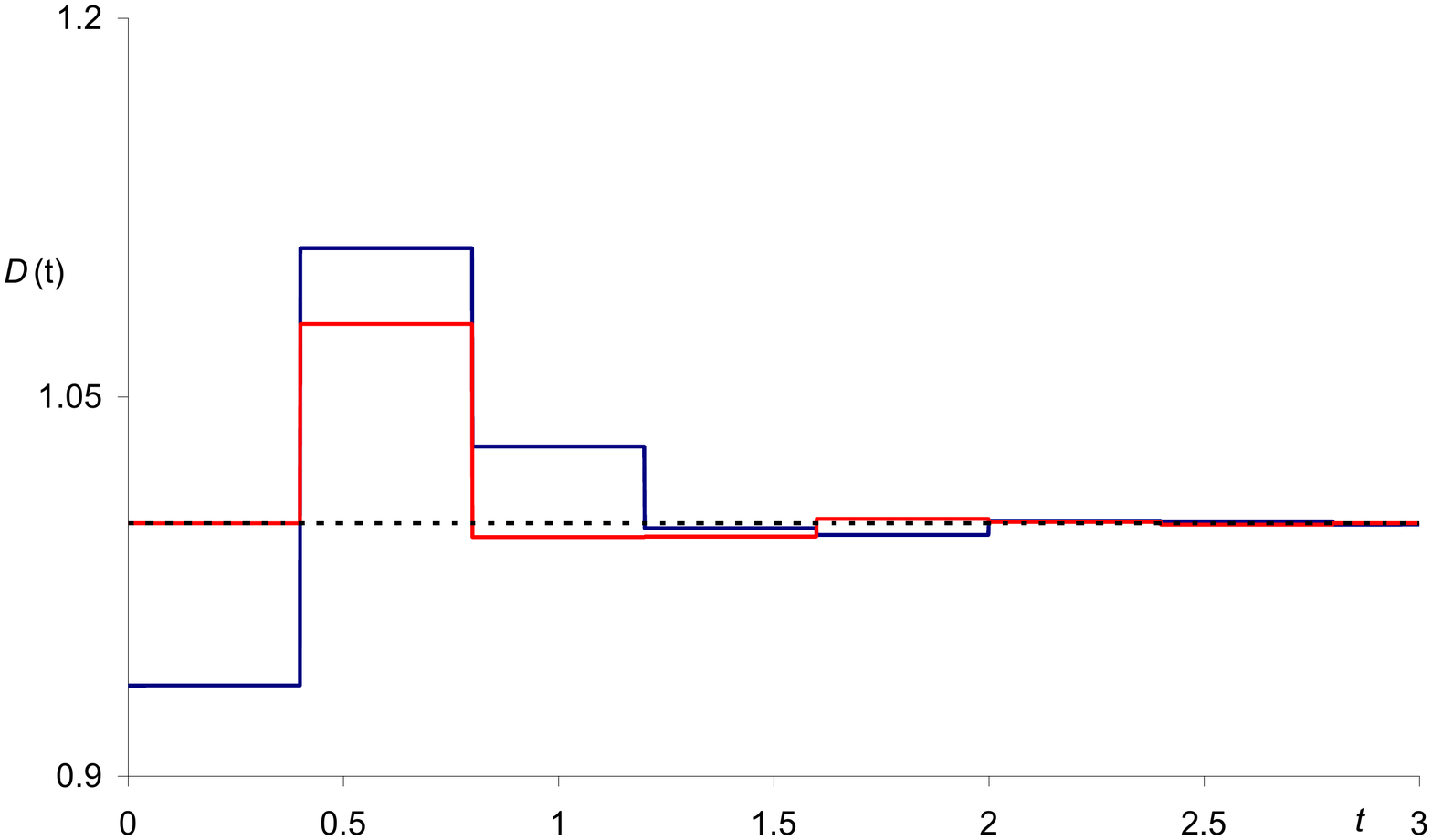}\\[-5.5em]
\caption{ First simulation. The red line is for the output feedback \rref{GrindEQ__7a_},
  the blue line is for the output feedback \rref{GrindEQ__8a_} and the black line is for the open-loop system with $D(t)\equiv 1$.}
\vspace{-.45em}
\end{figure}

In our second simulation,  we changed the parameter values to $b_{0} =1$, $b_{1} =0.7592106$, $c=4$, and $\theta =1$ and plotted the same values as before, in Figure 2. The responses for the output feedback law \rref{GrindEQ__7a_} and the output feedback law \rref{GrindEQ__8a_} are almost identical.

\begin{figure}[t]\vspace{-4.5em}
\hspace{.12in}\includegraphics[scale=0.5]{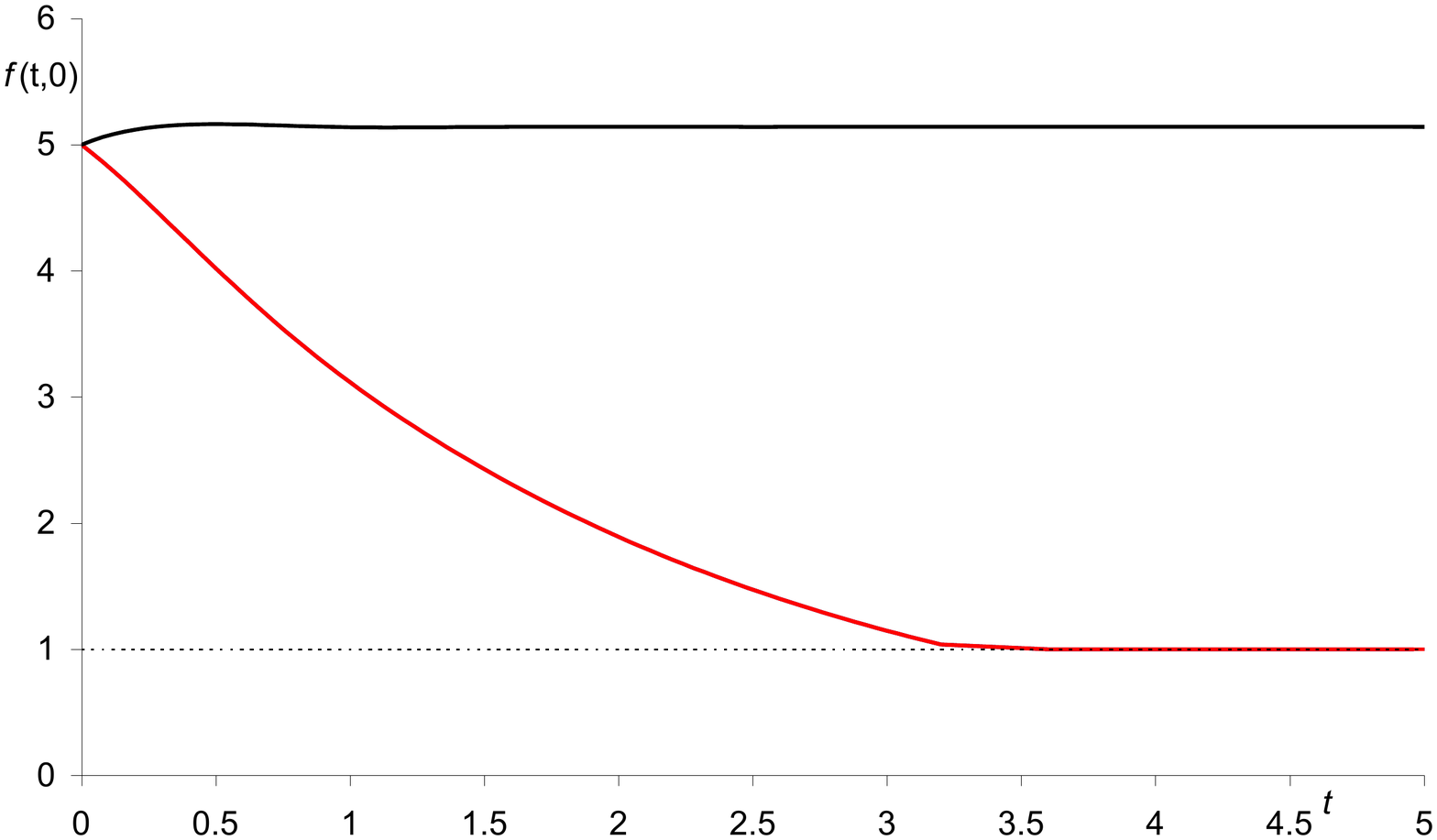}\\[-7em]
\includegraphics[scale=0.5]{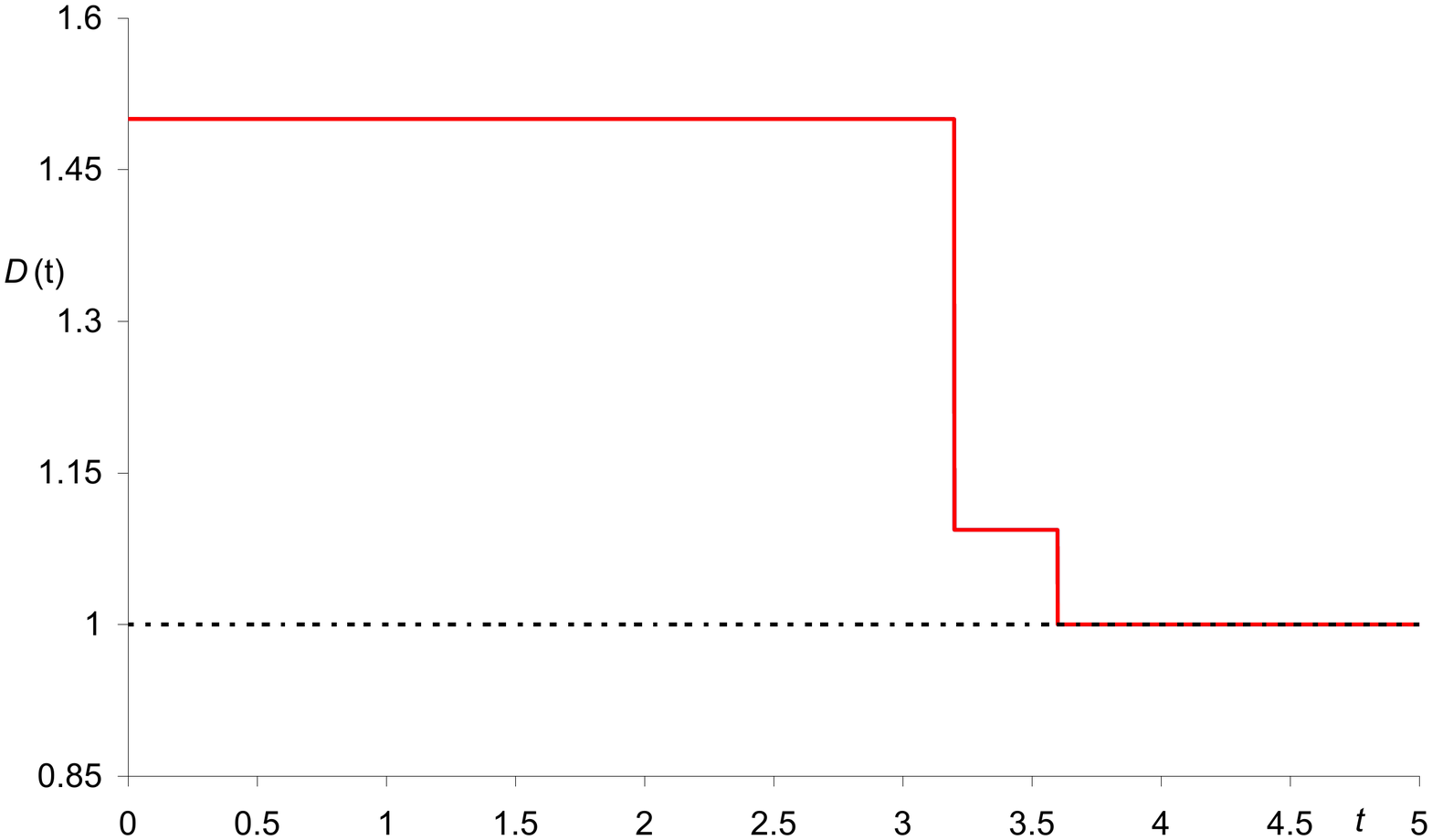}\\[-5.5em]
\caption{Second simulation. The red line is for the feedbacks \rref{GrindEQ__7a_} and  \rref{GrindEQ__8a_} and the black line is for the open-loop system with $D(t)\equiv 1$.}\vspace{-.15em}
\end{figure}

This second simulation was made  with an initial condition which is not close to the equilibrium profile (in the sense that it is an initial condition with very large initial population). The difference in the performance of the feedback controllers \rref{GrindEQ__7a_} and   \rref{GrindEQ__8a_} cannot be distinguished.

In our final simulation, we tested the robustness of the controller with respect to errors in the choice of $D^*$ being used in the controllers. We chose the   values $b_{0} =0.2$, $b_{1} =0.15184212$, $c=0.8$, and $\theta =1$, but instead of \rref{GrindEQ__7a_} and \rref{GrindEQ__8a_}, we applied the controllers which are
defined for all $t\in [iT,(i+1)T)$ and for  integers $i\ge 0$ by
\begin{equation}\label{7prime}D(t) =\max \left\{D_{\min } ,\min \left\{D_{\max } ,0.7+T^{-1} \ln \left(f(iT,0)/f^{*} (0)\right)\right\}\right\}  \end{equation}and
\begin{equation}\label{8prime} D(t) =\max \left\{D_{\min } ,\min \left\{D_{\max } ,0.7+T^{-1} \ln \left(y(iT)/y^{*} \right)\right\}\right\}.\end{equation}
In both cases, we obtained   $\lim_{t\to +\infty } f(t,0)=1.1275$ and 
 \begin{equation}
\lim_{t\to +\infty } f(t,0)=1.1275\; {\rm and}\; \lim_{t\to +\infty } D(t)=D^{*} =1.\end{equation}
Hence, a $-30\%$  error in $D^{*} $ gave a $+12.75\%$ steady-state deviation from the desired value of the newborn individuals. See Figure 3.

\begin{figure}[t]\vspace{-4.5em}
\includegraphics[scale=0.5]{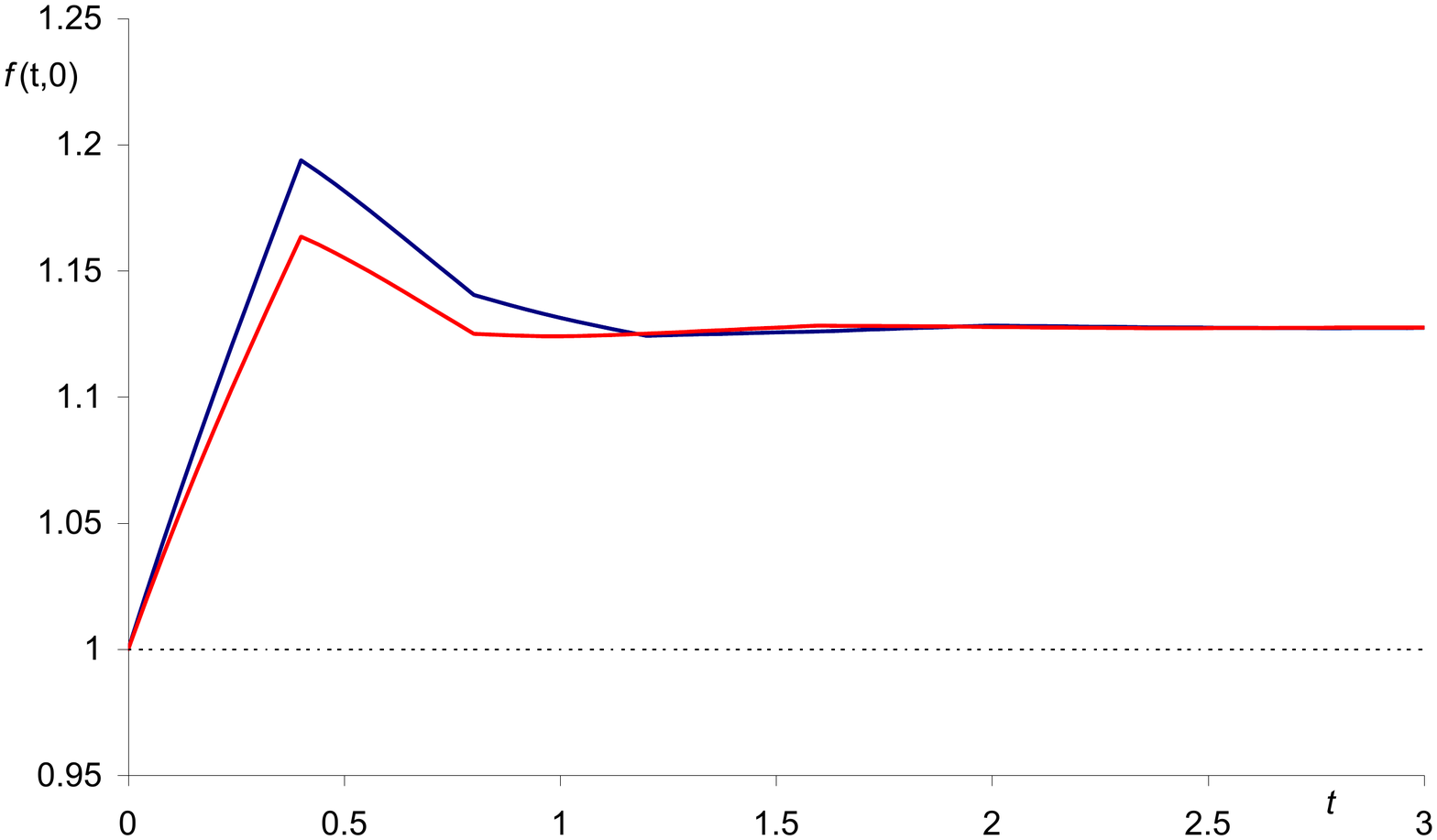}\\[-7em]
\includegraphics[scale=0.5]{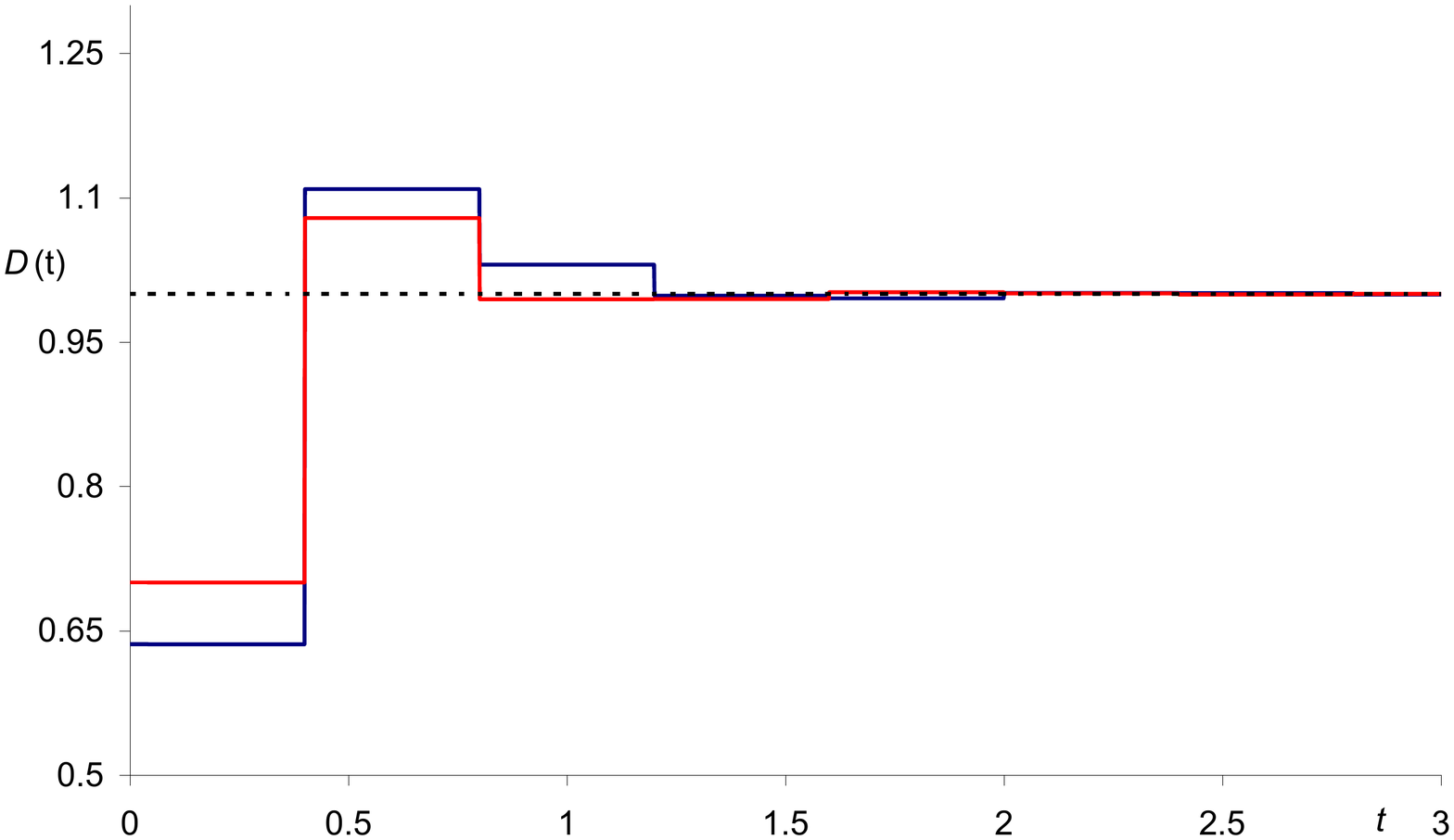}\\[-5.5em]
\caption{Third simulation. The red line is for the output feedback \rref{7prime}
  and the blue line is for the output feedback \rref{8prime}.}
\end{figure}

Notice that a constant error in $D^{*} $ is equivalent to an error in the set point since we have:
\begin{equation} \asa\begin{array}{l} {D(t) =\max \left\{D_{\min } ,\min \left\{D_{\max } ,0.7+T^{-1} \ln \left(f(iT,0)/f^{*} (0)\right)\right\}\right\}} \\ {=\max \left\{D_{\min } ,\min \left\{D_{\max } ,D^{*} +T^{-1} \ln \left(f(iT,0)/f^{*} (0)\right)-T^{-1} 0.12\right\}\right\}} \\ {=\max \left\{D_{\min } ,\min \left\{D_{\max } ,D^{*} +T^{-1} \ln \left(f(iT,0)/(1.1275f^{*} (0))\right)\right\}\right\}} \end{array}\end{equation}
  for the output feedback case \rref{GrindEQ__7a_} and
\begin{equation} \ashalf\begin{array}{l} {D(t) =\max \left\{D_{\min } ,\min \left\{D_{\max } ,0.7+T^{-1} \ln \left(y(iT)/y^{*} \right)\right\}\right\}} \\ {=\max \left\{D_{\min } ,\min \left\{D_{\max } ,D^{*} +T^{-1} \ln \left(y(iT)/y^{*} \right)-T^{-1} 0.12\right\}\right\}} \\ {=\max \left\{D_{\min } ,\min \left\{D_{\max } ,D^{*} +T^{-1} \ln \left(y(iT)/(1.1275y^{*} )\right)\right\}\right\}} \end{array}\end{equation}
 for the output feedback case \rref{GrindEQ__8a_}.  An interesting feature of the closed-loop system is that \begin{equation}\lim_{t\to +\infty } D(t)=D^{*} =1.\end{equation}
  \mm{$\lim_{t\to +\infty } D(t)=D^{*} =1$.}
  It may be worth considering an adaptive strategy for the elimination of errors in $D^{*} $ (i.e., a hybrid strategy that adapts the applied value of $D^{*} $). We leave the search for such a strategy for future work.

\section{Conclusions}
Chemostats play a vital role in  biotechnological applications, such as the production of insulin and in waste water treatment plants. Age-structured chemostats produce  challenging control problems for first-order hyperbolic PDEs that are beyond the scope of the existing controller methods for ODEs.
We studied the problem of stabilizing   an equilibrium age profile in an age-structured chemostat, using the  dilution rate as the control. We built a sampled-data dilution rate feedback  control law that ensures stability under arbitrary physically meaningful initial conditions and arbitrarily sparse sampling. Our control does not require measurement of the whole age profile, or exact model knowledge. The proposed feedback  also applies under arbitrary input constraints. The proof of our main result is based on  (a) the strong ergodic theorem and (b) our approach from \cite{KK14} for transforming a first-order hyperbolic PDE into an integral delay equation. Our simulations demonstrate the good performance of our controllers.
\newline\indent We hope to build on our research, in two ways. First,
since   the growth of the microorganism may sometimes depend on the concentration of a substrate, it would be useful to solve the stabilization problem for an enlarged system that has one PDE for the age distribution,  coupled with one ODE for the substrate (as proposed in \cite{TK06}, in the context of studying limit cycles with constant dilution rates instead of a control). Second, it would be useful to extend our work to cases where the   control is subject to uncertainties, and then seek generalizations of our exponential stability estimate such as input-to-state stability under input   constraints and sampling. Finally, we hope to cover state constrained problems, which   add the  requirement that the states must stay in prescribed subsets of the state space for all nonnegative times, in addition to the nonnegativity requirements on the physical quantities.

\renewcommand{\theequation}{A.\arabic{equation}}
\renewcommand{\thetheorem}{A.\arabic{theorem}}
\renewcommand{\thesection}{A.\arabic{section}}
  \setcounter{equation}{0}
\setcounter{section}{0}

\section*{Appendices}

\section{Proof of Lemma \ref{lm4}}\label{app1}
Local existence and uniqueness for every initial condition $x_{\scriptscriptstyle 0} \in L^{\infty } \left([-A,0);\mathbb R \right)$ is guaranteed by \cite[Theorem 2.1]{KK14}. We define the   functions  \mm{$V_*(t)=\sup_{-A\le a<0} x(t+a)$ and $W(t)=\inf_{-A\le a<0} x(t+a)$}
\begin{equation}
V_*(t)=\sup_{-A\le a<0} x(t+a)\; \; {\rm and}\; \; W(t)=\inf_{-A\le a<0} x(t+a)\end{equation}
for
  all $t\ge 0$ for which the solution of \rref{newide} exists.
  Let $q>0$ and $t\ge 0$ be sufficiently small such that the solution exists on $[t,t+q)$. From  equation \rref{newide}, we get
\[\asaa
\begin{array}{l} V_*(t+q)=\mathop{\sup }\limits_{-A\le a<0} x(t+q+a) =\mathop{\sup }\limits_{q-A\le s<q} x(t+s)\\=\max \left\{\mathop{\sup }\limits_{q-A\le s<0} x(t+s),\mathop{\sup }\limits_{0\le s<q} x(t+s)\right\} \\ \le \max \left\{V_*(t),\mathop{\sup }\limits_{0\le s<q} \int_{\scriptscriptstyle 0}^{\scriptscriptstyle A}G(a)x(t+s-a){\rm d}a \right\} \\   = \max \left\{V_*(t),\mathop{\sup }\limits_{0\le s<q} \left\{\int_{\Delta }^{\scriptscriptstyle A}G(a)x(t+s-a){\rm d}a+\int_{\scriptscriptstyle 0}^{\Delta }G(a)x(t+s-a){\rm d}a \right\}\right\} \\ \le \max \left\{V_*(t),\mathop{\sup }\limits_{0\le s<q} \! \left\{\mathop{\sup }\limits_{s-A\le l<s-\Delta } \! x(t+l) \int_{\Delta }^{\scriptscriptstyle A}G(a){\rm d}a  +\! \! \mathop{\sup }\limits_{s-\Delta \le l<s} \! x(t+l)\!  \int_{\scriptscriptstyle 0}^{\Delta }G(a){\rm d}a \! \right\}\! \right\} \\ \le \max \left\{V_*(t),\mathop{\sup }\limits_{-A\le l<q-\Delta } x(t+l)\int_{\Delta }^{\scriptscriptstyle A}G(a){\rm d}a +\mathop{\sup }\limits_{-\Delta \le l<q} x(t+l)\int_{\scriptscriptstyle 0}^{\Delta }G(a){\rm d}a \right\}. \end{array}\]
 Using the facts that
\begin{equation}\label{Elbig} \begin{array}{l}\int_{\scriptscriptstyle 0}^{\scriptscriptstyle A} G(a){\rm d}a =L \ge 1\; \; {\rm  and}\; \; \int_{\scriptscriptstyle 0}^{\Delta}G(a){\rm d}a <1 \end{array}\end{equation}
  and assuming that $q\le \min \{\Delta ,A-\Delta \}$, it follows that
\begin{equation} \label{GrindEQ_54a_}  \; \; \;  \; \; \;
V_*(t+q) \le  \max \left\{V_*(t),V_*(t)\left(L-\int_{\scriptscriptstyle 0}^{\Delta }G(a){\rm d}a \right)+V_*(t+q)\int_{\scriptscriptstyle 0}^{\Delta }G(a){\rm d}a \right\},
\end{equation}which gives $V_*(t+q)\le \max\{b V_*(t), V_*(t)\}$ with $b$ as in the statement of the lemma, by separately considering the two possible values for the maximum on the right side of \rref{GrindEQ_54a_}. (We need the maximum $\max\{b V_*(t), V_*(t)\}$ instead of just $V_*(t)$,  to allow the possibility that $V_*(t)$ is nonpositive.)

Similarly, we can use the decomposition
 \[\ashalf\begin{array}{l}W(t+q)=\mathop{\inf }\limits_{-A\le a<0} x(t+q+a)=\mathop{\inf }\limits_{q-A\le s<q} x(t+s)\\=\min \left\{\mathop{\inf }\limits_{q-A\le s<0}
x(t+s),\mathop{\inf }\limits_{0\le s<q} x(t+s)\right\},\end{array}\] the definition of $W(t)$, equation \rref{newide}, and the fact that $q\le \min \{\Delta ,A-\Delta \}$
to get \[\asb \begin{array}{l} W(t+q)\\\ge  \min \left\{W(t),\mathop{\inf }\limits_{0\le s<q} \! \int_{\scriptscriptstyle 0}^{\scriptscriptstyle A}G(a)x(t + s - a){\rm d}a
\right\} \\ =\min \left\{W(t),\mathop{\inf }\limits_{0\le s<q}  \left\{\int_{\Delta }^{\scriptscriptstyle A}G(a)x(t+s-a){\rm d}a+\int_{\scriptscriptstyle 0}^{\Delta
}G(a)x(t+s-a){\rm d}a \! \right\}\! \right\} \\ \ge \min \left\{W(t),\mathop{\inf }\limits_{0\le s<q}\!  \left\{\!  \mathop{\inf }\limits_{s-A\le l<s-\Delta } \! \! \! \! x(t+l) \int_{\Delta
}^{\scriptscriptstyle A}G(a){\rm d}a +\! \! \! \mathop{\inf }\limits_{s-\Delta \le l<s}\! \!  \! \! x(t+l) \int_{\scriptscriptstyle 0}^{\Delta }G(a){\rm d}a \! \right\}\! \right\} \\ \ge \min
\left\{W(t),\mathop{\inf }\limits_{-A\le l<q-\Delta } \! x(t+l)\int_{\Delta }^{\scriptscriptstyle A}G(a){\rm d}a + \mathop{\inf }\limits_{-\Delta \le l<q} x(t+l)\int_{\scriptscriptstyle
0}^{\Delta }G(a){\rm d}a \right\} \\ \ge \min \left\{W(t),W(t)\left(L-\int_{\scriptscriptstyle 0}^{\Delta }G(a){\rm d}a \right)+W(t+q)\int_{\scriptscriptstyle 0}^{\Delta }G(a){\rm
d}a \right\}, \end{array}\]
so \rref{Elbig} gives $W(t+q) \ge  \min\{b W(t), W(t)\}$.
 It follows from \rref{GrindEQ_54a_} that the solution of \rref{newide} is bounded  on $[t,t+q)$ when $q\le \min \{\Delta ,A-\Delta \}$. A standard contradiction argument in
conjunction with \cite[Theorem 2.1]{KK14} implies that the solution exists for all $t\ge 0$. Using induction and \rref{GrindEQ_54a_}, we can now show that \begin{equation} \min\{bW(0), W(0)\}\;
\le\; W(ih)\; \le\; V_*(ih)\; \le \; \max\{b^iV_*(0), V_*(0)\}
\label{GrindEQ_58_}\end{equation} for all integers $i\ge 0$, where $h=\min \{\Delta ,A-\Delta \}$. Inequality \rref{GrindEQ_14_} now follows from the definitions of $V_*$ and $W$ and
\rref{GrindEQ_54a_} and \rref{GrindEQ_58_}, by choosing $i$ such that $i\le t/h$. This proves the lemma.

\section{Proof of Claim \ref{claim1}}
\label{proofofclaim1}
We distinguish between the following cases.

  Case 1: $D_{\min } \le D^{*} +T^{-1} x_{i} \le D_{\max } $.
In this case, our choices of the $D_i$'s and $x_i$'s imply that   $D_{i} =D^{*} +T^{-1} x_{i} $. Using
our expressions \rref{GrindEQ_79_} for the output $y(t)$,  our choices of the $D_i$'s, the fact that $D_{i} =D^{*} +T^{-1} x_{i} $, and the fact $\phi (t)$ and $y(t)$ are continuous mappings, we get $x_{i+1} =u_{i} $, which directly implies \rref{GrindEQ_88_}.

 Case 2: $D_{\min } >D^{*} +T^{-1} x_{i} $.
In this case,  \rref{GrindEQ_86_} and  \rref{GrindEQ_87_} give $D_{i} =D_{\min } $.
Using the continuity of $\phi (t)$ and $y(t)$ and
setting $t=(i+1)T$ in   \rref{GrindEQ_79_}, and using definition \rref{GrindEQ_86_} and  the fact that $D_{i} =D_{\min } $, we get \[x_{i+1} =x_{i} +(D^{*} -D_{\min } )T+u_{i} .\] The inequality $D_{\min } >D^{*} +T^{-1} x_{i} $ implies that $0>(D^{*} -D_{\min } )T+x_{i} $ and $0>x_{i} $. Using the previous inequalities, the equality $x_{i+1} =x_{i} +(D^{*} -D_{\min } )T+u_{i} $, and our choice of $\delta$, we get
\[\ashalf\begin{array}{rcl} \left|x_{i+1} \right|&\le& \left|x_{i} +(D^{*} -D_{\min } )T\right|+\left|u_{i} \right| \; \,   =\; \,  -x_{i} -(D^{*} -D_{\min } )T+\left|u_{i} \right|\\ &\le &   |x_i|-2\delta+|u_i|, \end{array}\]
which again gives  \rref{GrindEQ_88_}.

  Case 3: $D^{*} +T^{-1} x_{i} >D_{\max } $. Then \rref{GrindEQ_86_} and  \rref{GrindEQ_87_}
   give $D_{i} =D_{\max } $. Since $\phi$ and $y$ are continuous, we can set   $t=(i+1)T$ in   \rref{GrindEQ_79_},
  and use   \rref{GrindEQ_86_} and  the fact that $D_{i} =D_{\min } $ to get $x_{i+1} =x_{i} -(D_{\max } -D^{*} )T+u_{i} $. The inequality $D^{*} +T^{-1} x_{i} >D_{\max } $ implies that $x_{i} -(D_{\max } -D^{*} )T>0$ and $0<x_{i} $. Then  the equality $x_{i+1} =x_{i} -(D_{\max } -D^{*} )T+u_{i} $ give
\[\ashalf\begin{array}{rcl}  \left|x_{i+1} \right|&\le& \left|x_{i} -(D_{\max } -D^{*} )T\right|+\left|u_{i} \right| \;  \,  =\; \, x_{i} -(D_{\max } -D^{*} )T+\left|u_{i} \right| \\& \le & |x_i| -2\delta+|u_i|,
\end{array}\]
so \rref{GrindEQ_88_} holds again. This proves Claim \ref{claim1}.

\section{Proof of Claim \ref{claim2}}
\label{proofofclaim2}
 The proof of  \rref{GrindEQ_89_} is  by induction. First notice that both inequalities in \rref{GrindEQ_89_} hold for $i=0$. Next assume that  \rref{GrindEQ_89_} hold for certain integer $i\ge 0$. We consider three cases.

Case 1: $D_{\min } \le D^{*} +T^{-1} x_{i} \le D_{\max } $. In this case, our treatment of Case 1 in our proof of Claim \ref{claim1} gives $x_{i+1} =u_{i} $. Consequently,  our definition \rref{GrindEQ_87_} gives
\[\asa\begin{array}{l} x_{i+1} =u_{i} =\ln \left(\frac{P(f_{0} )+\int _{0}^{A}g(a)\phi ((i+1)T-a){\rm d}a }{P(f_{0} )+\int _{0}^{A}g(a)\phi (iT-a){\rm d}a } \right)\\\le \max_{k=0,...,i+1} \left(\ln \left(\frac{P(f_{0} )+\int _{0}^{A}g(a)\phi ((i+1)T-a){\rm d}a }{P(f_{0} )+\int _{0}^{A}g(a)\phi (kT-a){\rm d}a } \right)\right) \\ \le \max \left\{0,x_{0} -(i+1)\left(D_{\max } -D^{*} \right)T\right\}\\+\max_{k=0,...,i+1} \left(\ln \left(\frac{P(f_{0} )+\int _{0}^{A}g(a)\phi ((i+1)T-a){\rm d}a }{P(f_{0} )+\int _{0}^{A}g(a)\phi (kT-a){\rm d}a } \right)\right), \end{array}\]which
implies the second inequality in \rref{GrindEQ_89_} with $i+1$ in place of $i\ge 0$. Similarly, we obtain the first inequality in \rref{GrindEQ_89_} with $i+1$ in place of $i\ge 0$.

Case 2: $D_{\min } >D^{*} +T^{-1} x_{i} $. Arguing as in our treatment of Case 2 in our proof of Claim \ref{claim1}, we get $x_{i+1} =x_{i} +(D^{*} -D_{\min } )T+u_{i} $. Hence,  \rref{GrindEQ_87_} and \rref{GrindEQ_89_} and the fact that $D^*\ge D_{\rm min}$ give
\[\asa\begin{array}{rcl} x_{i+1} &=&x_{i} +\left(D^{*} -D_{\min } \right)T+u_{i} \\&\ge& \min \left\{0,x_{0} +i\left(D^{*} -D_{\min } \right)T\right\}+\left(D^{*} -D_{\min } \right)T \\ &&+\min_{k=0,...,i} \left(\ln \left(\frac{P(f_{0} )+\int _{0}^{A}g(a)\phi (iT-a){\rm d}a }{P(f_{0} )+\int _{0}^{A}g(a)\phi (kT-a){\rm d}a } \right)\right)\\&&+\, \ln \left(\frac{P(f_{0} )+\int _{0}^{A}g(a)\phi ((i+1)T-a){\rm d}a }{P(f_{0} )+\int _{0}^{A}g(a)\phi (iT-a){\rm d}a } \right) \\ &\ge& \min \left\{\left(D^{*} -D_{\min } \right)T,x_{0} +(i+1)\left(D^{*} -D_{\min } \right)T\right\}\\&&+\; \min_{k=0,...,i} \left(\ln \left(\frac{P(f_{0} )+\int _{0}^{A}g(a)\phi ((i+1)T-a){\rm d}a }{P(f_{0} )+\int _{0}^{A}g(a)\phi (kT-a){\rm d}a } \right)\right) \\ &\ge& \min \left\{0,x_{0} +(i+1)\left(D^{*} -D_{\min } \right)T\right\}\\&&+\, \min_{k=0,...,i+1} \left(\ln \left(\frac{P(f_{0} )+\int _{0}^{A}g(a)\phi ((i+1)T-a){\rm d}a }{P(f_{0} )+\int _{0}^{A}g(a)\phi (kT-a){\rm d}a } \right)\right) \end{array}\]
which is the first inequality \rref{GrindEQ_89_} with $i+1$ in place of $i\ge 0$. Furthermore,  $0>(D^{*} -D_{\min } )T+x_{i} $. Combining the previous inequality with definition \rref{GrindEQ_87_} and the fact that $x_{i+1} =x_{i} +(D^{*} -D_{\min } )T+u_{i} $, we get
\[\asa\begin{array}{rcl} x_{i+1} &=&x_{i} +\left(D^{*} -D_{\min } \right)T+u_{i} \le u_{i}\\& \le &  \max_{k=0,...,i+1} \left(\ln \left(\frac{P(f_{0} )+\int _{0}^{A}g(a)\phi ((i+1)T-a){\rm d}a }{P(f_{0} )+\int _{0}^{A}g(a)\phi (kT-a){\rm d}a } \right)\right) \\ &\le& \max \left\{0,x_{0} -(i+1)\left(D_{\max } -D^{*} \right)T\right\}\\&&+\,  \max_{k=0,...,i+1} \left(\ln \left(\frac{P(f_{0} )+\int _{0}^{A}g(a)\phi ((i+1)T-a){\rm d}a }{P(f_{0} )+\int _{0}^{A}g(a)\phi (kT-a){\rm d}a } \right)\right) \end{array}\]
 which is the second inequality \rref{GrindEQ_89_} with $i+1$ in place of $i\ge 0$.

Case 3: $D^{*} +T^{-1} x_{i} >D_{\max } $. Arguing as in Case 3 in the proof of Claim \ref{claim1} gives $x_{i+1} =x_{i} -(D_{\max } -D^{*} )T+u_{i} $. Combining the previous equality with definition \rref{GrindEQ_87_} and inequality \rref{GrindEQ_89_}, we get
\[\asa\begin{array}{rcl} x_{i+1} &=&x_{i} -\left(D_{\max } -D^{*} \right)T+u_{i} \\&\le& \max \left\{0,x_{0} -i\left(D_{\max } -D^{*} \right)T\right\}-\left(D_{\max } -D^{*} \right)T \\ &&+\max_{k=0,...,i} \left(\ln \left(\frac{P(f_{0} )+\int _{0}^{A}g(a)\phi (iT-a){\rm d}a }{P(f_{0} )+\int _{0}^{A}g(a)\phi (kT-a){\rm d}a } \right)\right)\\&&+\, \ln \left(\frac{P(f_{0} )+\int _{0}^{A}g(a)\phi ((i+1)T-a){\rm d}a }{P(f_{0} )+\int _{0}^{A}g(a)\phi (iT-a){\rm d}a } \right) \\ &\le& \max \left\{-\left(D_{\max } -D^{*} \right)T,x_{0} -(i+1)\left(D_{\max } -D^{*} \right)T\right\}\\&&+\, \max_{k=0,...,i} \left(\ln \left(\frac{P(f_{0} )+\int _{0}^{A}g(a)\phi ((i+1)T-a){\rm d}a }{P(f_{0} )+\int _{0}^{A}g(a)\phi (kT-a){\rm d}a } \right)\right) \\ &\le& \max \left\{0,x_{0} -(i+1)\left(D_{\max } -D^{*} \right)T\right\}\\&&+\, \max_{k=0,...,i+1} \left(\ln \left(\frac{P(f_{0} )+\int _{0}^{A}g(a)\phi ((i+1)T-a){\rm d}a }{P(f_{0} )+\int _{0}^{A}g(a)\phi (kT-a){\rm d}a } \right)\right), \end{array}\]
which is the second inequality \rref{GrindEQ_89_} with $i+1$ in place of $i\ge 0$. Furthermore,   $x_{i} -(D_{\max } -D^{*} )T>0$. Combining the previous inequality with definition \rref{GrindEQ_87_} and the fact that $x_{i+1} =x_{i} -(D_{\max } -D^{*} )T+u_{i} $, we get
\[\asa\begin{array}{rcl} x_{i+1} &=&x_{i} -\left(D_{\max } -D^{*} \right)T+u_{i} \\& \ge&  u_{i} \;  \, \ge\;  \,  \min_{k=0,...,i+1} \left(\ln \left(\frac{P(f_{0} )+\int _{0}^{A}g(a)\phi ((i+1)T-a){\rm d}a }{P(f_{0} )+\int _{0}^{A}g(a)\phi (kT-a){\rm d}a } \right)\right) \\ &\ge& \min \left\{0,x_{0} +(i+1)\left(D^{*} -D_{\min } \right)T\right\}\\&&+\, \min_{k=0,...,i+1} \left(\ln \left(\frac{P(f_{0} )+\int _{0}^{A}g(a)\phi ((i+1)T-a){\rm d}a }{P(f_{0} )+\int _{0}^{A}g(a)\phi (kT-a){\rm d}a } \right)\right), \end{array}\]
  which is the first inequality \rref{GrindEQ_89_} with $i+1$ in place of $i\ge 0$. This proves Claim \ref{claim2}.

\section{Proof of Claim \ref{claim3}}\label{proofofclaim3}
Our expressions  \rref{GrindEQ_79_} for the output give
\begin{equation} \label{GrindEQ_91_}\begin{array}{l}
x(t)=x_{i} -(D_{i} -D^{*} )(t-iT)+\ln \left(\frac{P(f_{0} )+\int _{0}^{A}g(a)\phi (t-a){\rm d}a }{P(f_{0} )+\int _{0}^{A}g(a)\phi (iT-a){\rm d}a } \right)\end{array}
\end{equation}
for all $t\ge 0$,
  where $i=[t/T]$. We again consider three cases.

\noindent Case 1: $D_{\min } \le D^{*} +T^{-1} x_{i} \le D_{\max } $. In this case,
 \rref{GrindEQ_86_} and  \rref{GrindEQ_87_} imply that $D_{i} =D^{*} +T^{-1} x_{i} $. Hence, \rref{GrindEQ_91_}gives
\[\begin{array}{l}x(t)=\left(1-(t-iT)T^{-1} \right)x_{i} +\ln \left(\frac{P(f_{0} )+\int _{0}^{A}g(a)\phi (t-a){\rm d}a }{P(f_{0} )+\int _{0}^{A}g(a)\phi (iT-a){\rm d}a } \right)\end{array}\]
for all $t\ge 0$.  The above equality in conjunction with the fact that $0\le 1-(t-iT)T^{-1} \le 1$ and inequality \rref{GrindEQ_89_} gives estimates \rref{90a}-\rref{90b}.

 Case 2: $D_{\min } >D^{*} +T^{-1} x_{i} $. Now our definitions \rref{GrindEQ_86_} and  \rref{GrindEQ_87_} give   $D_{i} =D_{\min } $. The inequality $D_{\min } >D^{*} +T^{-1} x_{i} $ implies that $0>(D^{*} -D_{\min } )T+x_{i} $ and $0>x_{i} $, which combined with \rref{GrindEQ_91_} give:
\[\begin{array}{l}x_{i} +\ln \left(\frac{P(f_{0} )+\int _{0}^{A}g(a)\phi (t-a){\rm d}a }{P(f_{0} )+\int _{0}^{A}g(a)\phi (iT-a){\rm d}a } \right)\le x(t)\le \ln \left(\frac{P(f_{0} )+\int _{0}^{A}g(a)\phi (t-a){\rm d}a }{P(f_{0} )+\int _{0}^{A}g(a)\phi (iT-a){\rm d}a } \right).\end{array}\]
 The above inequality in conjunction with inequalities \rref{GrindEQ_89_} gives \rref{90a}-\rref{90b}.

  Case 3: $D^{*} +T^{-1} x_{i} >D_{\max } $. Definitions \rref{GrindEQ_86_}-\rref{GrindEQ_87_} imply that $D_{i} =D_{\max } $. The inequality $D^{*} +T^{-1} x_{i} >D_{\max } $ implies that $x_{i} -(D_{\max } -D^{*} )T>0$ and $0<x_{i} $, which combined with \rref{GrindEQ_91_} gives
\[\begin{array}{l}\ln \left(\frac{P(f_{0} )+\int _{0}^{A}g(a)\phi (t-a){\rm d}a }{P(f_{0} )+\int _{0}^{A}g(a)\phi (iT-a){\rm d}a } \right)\le x(t)\le x_{i} +\ln \left(\frac{P(f_{0} )+\int _{0}^{A}g(a)\phi (t-a){\rm d}a }{P(f_{0} )+\int _{0}^{A}g(a)\phi (iT-a){\rm d}a } \right),\end{array}\]
which we can combine with \rref{GrindEQ_89_} to get  \rref{90a}-\rref{90b}. This proves Claim \ref{claim3}.

\section{Proof of Claim \ref{claim4}}\label{proofofclaim4}
Since $\left|\ln (x)\right|=\ln \left(\max \{x,x^{-1} \}\right)$ for all $x>0$, we get
\[\asa\begin{array}{l}\left|\ln \left(\frac{P(f_{0} )+\int _{0}^{A}g(a)\phi (t-a){\rm d}a }{P(f_{0} )+\int _{0}^{A}g(a)\phi (iT-a){\rm d}a } \right)\right|=\\\ln \left(\max \left\{\frac{P(f_{0} )+\int _{0}^{A}g(a)\phi (t-a){\rm d}a }{P(f_{0} )+\int _{0}^{A}g(a)\phi (iT-a){\rm d}a } \, ,\, \frac{P(f_{0} )+\int _{0}^{A}g(a)\phi (iT-a){\rm d}a }{P(f_{0} )+\int _{0}^{A}g(a)\phi (t-a){\rm d}a } \right\}\right)\end{array}\]
and \begin{equation}\; \; \; \; \begin{array}{l}\left|u_{i} \right|=\ln \left(\max \left\{\frac{P(f_{0} )+\int _{0}^{A}g(a)\phi ((i+1)T-a){\rm d}a }{P(f_{0} )+\int _{0}^{A}g(a)\phi (iT-a){\rm d}a } \, ,\, \frac{P(f_{0} )+\int _{0}^{A}g(a)\phi (iT-a){\rm d}a }{P(f_{0} )+\int _{0}^{A}g(a)\phi ((i+1)T-a){\rm d}a } \right\}\right).  \end{array}   \label{GrindEQ_94_}\end{equation}
  On the other hand, we can use  \rref{GrindEQ_93_} to get the following for all $i\ge j$:
\begin{equation} \label{GrindEQ_95_}\begin{array}{l}
\frac{P(f_{0} )+h}{P(f_{0} )-h} \le \exp (\delta ),\end{array}
\end{equation}
  where $h=K^{*} \left\| f_{0} \right\| _{1} \exp \left(-\varepsilon (iT-A)\right)$. Using
  our bound \rref{GrindEQ_35_} on $\phi$
   and the fact $j\ge [A/T]+1$ (which implies that $jT\ge A$, i.e., $iT-a\ge 0$ for all $a\in [0,A]$), we get
\[\asa\begin{array}{l} \frac{P(f_{0} )-K^{*} \left\| f_{0} \right\| _{1} \exp \left(-\varepsilon ((i+1)T-A)\right)\int _{0}^{A}g(a){\rm d}a }{P(f_{0} )+K^{*} \left\| f_{0} \right\| _{1} \exp \left(-\varepsilon (iT-A)\right)\int _{0}^{A}g(a){\rm d}a } \le    \frac{P(f_{0} )+\int _{0}^{A}g(a)\phi ((i+1)T-a){\rm d}a }{P(f_{0} )+\int _{0}^{A}g(a)\phi (iT-a){\rm d}a }\\ \le\,  \frac{P(f_{0} )+K^{*} \left\| f_{0} \right\| _{1} \exp \left(-\varepsilon ((i+1)T-A)\right)\int _{0}^{A}g(a){\rm d}a }{P(f_{0} )-K^{*} \left\| f_{0} \right\| _{1} \exp \left(-\varepsilon (iT-A)\right)\int _{0}^{A}g(a){\rm d}a }  \end{array}\]
and \[\begin{array}{l} \frac{P(f_{0} )-K^{*} \left\| f_{0} \right\| _{1} \exp \left(-\varepsilon (t-A)\right)\int _{0}^{A}g(a){\rm d}a }{P(f_{0} )+K^{*} \left\| f_{0} \right\| _{1} \exp \left(-\varepsilon (iT-A)\right)\int _{0}^{A}g(a){\rm d}a }\\ \le \frac{P(f_{0} )+\int _{0}^{A}g(a)\phi (t-a){\rm d}a }{P(f_{0} )+\int _{0}^{A}g(a)\phi (iT-a){\rm d}a } \le \frac{P(f_{0} )+K^{*} \left\| f_{0} \right\| _{1} \exp \left(-\varepsilon (t-A)\right)\int _{0}^{A}g(a){\rm d}a }{P(f_{0} )-K^{*} \left\| f_{0} \right\| _{1} \exp \left(-\varepsilon (iT-A)\right)\int _{0}^{A}g(a){\rm d}a }  \end{array}\] for all $t\ge iT$.

Our formulas \rref{GrindEQ_81_} and \rref{GrindEQ_84_} for $\beta$ and $g$ imply that
 $\int _{0}^{A}g(a){\rm d}a =1$. Hence, the preceding inequalities give
\[\begin{array}{l}\frac{P(f_{0} )-h}{P(f_{0} )+h} \le \frac{P(f_{0} )+\int _{0}^{A}g(a)\phi ((i+1)T-a){\rm d}a }{P(f_{0} )+\int _{0}^{A}g(a)\phi (iT-a){\rm d}a } \le \frac{P(f_{0} )+h}{P(f_{0} )-h} \end{array}\] and
\begin{equation}\begin{array}{l}\frac{P(f_{0} )-h}{P(f_{0} )+h} \le \frac{P(f_{0} )+\int _{0}^{A}g(a)\phi (t-a){\rm d}a }{P(f_{0} )+\int _{0}^{A}g(a)\phi (iT-a){\rm d}a } \le \frac{P(f_{0} )+h}{P(f_{0} )-h} \; \; {\rm  for\ all}\; t\ge iT   \; ,           \end{array}\label{GrindEQ_96_}\end{equation}
  where $h=K^{*} \left\| f_{0} \right\| _{1} \exp \left(-\varepsilon (iT-A)\right)$. Combining \rref{GrindEQ_94_} and  \rref{GrindEQ_96_}, we get:
\begin{equation}\begin{array}{l}\left|u_{i} \right|\le \ln \left(\frac{P(f_{0} )+h}{P(f_{0} )-h} \right)\; \; {\rm and}\\ \left|\ln \left(\frac{P(f_{0} )+\int _{0}^{A}g(a)\phi (t-a){\rm d}a }{P(f_{0} )+\int _{0}^{A}g(a)\phi (iT-a){\rm d}a } \right)\right|\le \ln \left(\frac{P(f_{0} )+h}{P(f_{0} )-h} \right)\; \; {\rm  for\ all}\; t\ge iT \; . \end{array}\label{GrindEQ_97_}\end{equation}
 Using \rref{GrindEQ_95_} and \rref{GrindEQ_97_} we obtain the desired inequality $\left|u_{i} \right|\le \delta $ for all $i\ge j$. Also, our assumptions on $k$ and $f_0$ ensure that $P(f_0)>0$.

  Next, using \rref{GrindEQ_93_} and  \rref{GrindEQ_97_}, writing
  \begin{equation}\begin{array}{l}\frac{P(f_0)+h}{P(f_0)-h}=1+\frac{2h}{h-P(f_0)},\end{array}\end{equation}
  and using the  inequality $\ln (1+x)\le x$ for all $x\ge 0$, we obtain:
\[\left|u_{i} \right|\le \frac{2h}{P(f_{0} )-h} \le \frac{2h}{P(f_{0} )-\frac{\exp (\delta )-1}{\exp (\delta )+1} P(f_{0} )} =\frac{h\left(\exp (\delta )+1\right)}{P(f_{0} )}\]
and
\begin{equation}\begin{array}{l}\left|\ln \left(\frac{P(f_{0} )+\int _{0}^{A}g(a)\phi (t-a){\rm d}a }{P(f_{0} )+\int _{0}^{A}g(a)\phi (iT-a){\rm d}a } \right)\right|\le \frac{h\left(\exp (\delta )+1\right)}{P(f_{0} )} \; \; {\rm  for\ all}\; t\ge iT .                    \end{array} \label{GrindEQ_98_}\end{equation}
  Since $h=K^{*} \left\| f_{0} \right\| _{1} \exp \left(-\varepsilon (iT-A)\right)$,  we obtain
  \begin{equation}
  \begin{array}{l}\left|u_{i} \right|\le \frac{K^{*} \left\| f_{0} \right\| _{1} (\exp (\delta )+1)\exp (\varepsilon A)}{P(f_{0} )} \exp (-\varepsilon \, iT)\; \;  {\rm  for\  all}\;  i\ge j\; \; {\rm and}\end{array}\end{equation}
  \begin{equation}\begin{array}{l}
   \left|\ln \left(\frac{P(f_{0} )+\int _{0}^{A}g(a)\phi (t-a){\rm d}a }{P(f_{0} )+\int _{0}^{A}g(a)\phi (iT-a){\rm d}a } \right)\right|\le \frac{K^{*} \left\| f_{0} \right\| _{1} \left(\exp (\delta )+1\right)\exp (\varepsilon A)}{P(f_{0} )} \exp \left(-\varepsilon \, iT\right)\end{array}\end{equation}
  for all  $i\ge j$ and  $t\ge iT$.
This completes the proof of Claim \ref{claim4}.

\end{document}